\documentclass[10pt]{amsart}

\usepackage{amsmath}
\usepackage{amssymb}
\usepackage{amsfonts}
\usepackage{geometry}
\usepackage{url}

\usepackage{amsmath}
\usepackage{amssymb}
\usepackage{amsfonts}
\usepackage{geometry}
\usepackage{url}
\usepackage{setspace}
\usepackage{amsthm}
\usepackage[all,cmtip]{xy}
\usepackage{amsmath,amscd}

\usepackage{amsthm}
\newtheorem{theorem}{Theorem}[section]
\newtheorem{lemma}[theorem]{Lemma}
\newtheorem{proposition}[theorem]{Proposition}

\theoremstyle{definition}
\newtheorem{definition}[theorem]{Definition}

\newtheorem{remark}[theorem]{Remark}

\newcommand{\Hom}{\mathrm{Hom}}
\newcommand{\Sym}{\mathrm{S}}
\newcommand{\Res}{\mathrm{Res}}
\newcommand{\Ind}{\mathrm{Ind}}

\begin{document}

\title{Irreducible representations of the rational Cherednik algebra associated to the Coxeter group $H_3$}

\author{Martina Balagovi\' c \and Arjun Puranik}
\address{Department of Mathematics,  Massachusetts Institute
of Technology, Cambridge, MA 02139, USA }
\email{martinab@math.mit.edu, arap1234@gmail.com}
\maketitle

\maketitle

\begin{abstract} This paper describes irreducible representations in category $\mathcal{O}$ of the rational Cherednik algebra $H_{c}(H_3,\mathfrak{h})$ associated to the exceptional Coxeter group $H_3$ and any complex parameter $c$. We compute the characters of all these representations explicitly. As a consequence, we classify all the finite dimensional irreducible representations of  $H_{c}(H_3,\mathfrak{h})$. 
\end{abstract}

\section{ Introduction}

A rational Cherednik algebra $H_c(W,\mathfrak{h})$ is a certain associative algebra defined by a finite Coxeter group $W$, its complexified reflection representation $\mathfrak{h}$ and a certain parameter $c$. In case $W$ is a Weyl group, these algebras are rational degenerations of double affine Hecke algebras, which were defined by I. Cherednik \cite{Ch} and used to prove Macdonald conjectures. They can also be thought of in relation to completely integrable systems as algebras encoding the structure of Dunkl operators \cite{D}, \cite{DO}, and Calogero-Moser systems \cite{E},  or as a special case of symplectic reflection algebras of Etingof and Ginzburg \cite{EG}. These algebras and their representation theory have been intensively studied in the last fifteen years.

For  $H_c(W,\mathfrak{h})$ one can define category $\mathcal{O}$ (sometimes denoted $\mathcal{O}_{c}(W,\mathfrak{h})$ or $\mathcal{O}_{c}$, when there is need to explicitly mention $W,\mathfrak{h},c$), which has some simiSSlarities to the very well understood category $\mathcal{O}$ from Lie theory \cite{GGOR}. For example, one can define standard (Verma) modules $M_{c}(\tau)$ (which depend on parameter $\tau$, an irreducible representation of $W$, which appears instead of the lowest weight vector). Every $M_{c}(\tau)$ has a unique irreducible quotient $L_{c}(\tau)$. These are the only irreducible modules in category $\mathcal{O}$. There is a contravariant form $B$ analogous to  Shapovalov form for Lie algebra representations, defined on $M_{c}(\tau)$ and nondegenerate on its quotient $L_{c}(\tau)$. Characters can be defined which determine the irreducible modules completely, and are linearly independent. 

An obviously interesting question that appears is to describe all irreducible objects in category $\mathcal{O}_{c}(W,\mathfrak{h})$, meaning all $L_{c}(\tau)$. This can be done for example by writing their characters or by giving a description in the Grothendieck group in terms of $M_{c}(\tau)$, in a way analogous to the Weyl character formula. For rational Cherednik algebras such a description doesn't yet exist; it is not even known in general which of the  $L_{c}(\tau)$ are finite dimensional. 

Many partial results exist. For $(W,\mathfrak{h})$ of type $A$, \cite{BEG1} calculates the character formulas for all the finite dimensional $L_{c}(\tau)$. Also for type $A$, \cite{R1} calculates all the characters for $c$ not a half integer, and conjectures that the analogous formulas hold for $c$ a half integer. For dihedral groups, \cite{Chm} computes the characters of irreducible modules in category $\mathcal{O}$. The paper  \cite{VV} answers the question of when is the representation $L_{c}(\tau)$ finite dimensional for $W$ a Weyl group, $c$ a constant, and $\tau$ a trivial representation of $W$. A generalization of this is a recent result by Etingof (see \cite{E2}), which gives an answer for $W$ any finite Coxeter, trivial $\tau$, and any value of the parameter $c$. 

This paper considers the algebra $H_c(W,\mathfrak{h})$ for the case when $W$ is the exceptional Coxeter group $H_3$. In that case the parameter $c$ is just a complex number. We calculate the characters of all irreducible representations in category $\mathcal{O}$, i.e. $L_{c}(\tau)$ for all values of $c$ and $\tau$, and consequently classify the finite dimensional ones. The main result is in Theorem \ref{main}, which is proved, case by case, by Theorems \ref{equi1}, \ref{equi2}, \ref{1/10}, \ref{1/6}, \ref{1/5}, \ref{1/3}, \ref{1/2} and \ref{3/2}.

Calculating the characters provided us with an example $c=1/2$ of an \textit{aspherical} value for $c$; meaning such that there exists a module in category $\mathcal{O}_{c}$ that has trivial $W$-invariants (namely, $L_{1/2}(\widetilde{\mathbf{3}}_{-})$). It was conjectured (\cite{BE})  that such values of $c$ can only be in $\left( -1,0 \right)$. $c=1/2$ for $W=H_3$ was the first counterexample; counterexamples in type $B$ have been found by P. Etingof shortly afterwards. 

Our method is to first use results by Bezrukavnikov, Etingof and Rouquier to reduce the set of pairs $(c,\tau)$ for which we need to calculate the characters to a small finite set. Namely, category $\mathcal{O}_c(W,\mathfrak{h})$ is semisimple unless $c$ is a \textit{singular} parameter for $W$. It is known that if $c$ is a constant (as it is in our case), it is singular if and only if it  is a rational number whose denominator divides a degree of a basic invariant of $W$ (see \cite{GGOR} and \cite{BE}). Then, we use an equivalence of categories of representations of $H_c(W,\mathfrak{h})$ and $H_{-c}(W,\mathfrak{h})$ to reduce to the case of positive $c$ ($c=0$ is trivial).  Next, there are equivalences of categories between category $\mathcal{O}_{1/d}(W,\mathfrak{h})$ and category $\mathcal{O}_{r/d}(W,\mathfrak{h})$, in the case $d\ne 2$ (\cite{R2}), and finally between category $\mathcal{O}_{c}(W,\mathfrak{h})$ and category $\mathcal{O}_{c+1}(W,\mathfrak{h})$ for $c>>0$ \cite{BEG1}. It is known how these functors act on the standard and irreducible modules, and consequently how the characters transform under them. All this allows us to reduce the possible values of $c$ that we need to consider to a very small set; $c\in\lbrace 1/10,1/6,1/5,1/3,1/2,3/2\rbrace$.

So, we are left with the task of computing $L_c(\tau)$ for a small finite number of $c\in \mathbb{C}$ and all $\tau \in Irrep(W)$. We do this by more or less elementary methods, expressing the characters of $L_{c}(\tau)$ in terms of characters of various $M_{c}(\sigma)$. There is a copy of $\mathfrak{sl}_2$ in $H_c(W,\mathfrak{h})$, and it is such that its semisimple element acts as a grading element on both  $H_c(W,\mathfrak{h})$ and all its representations in category $\mathcal{O}$. There is also a copy of the group algebra $\mathbb{C}W$ in $H_c(W,\mathfrak{h})$, and it commutes with the action of the grading element. Thus any representation in category $\mathcal{O}$ has a grading with graded pieces being finite dimensional representations of $W$. We use simple observations of this type, along with computations in representation theory of $\mathfrak{sl}_2$ and $H_3$, to narrow down the options for possible coefficients in character formulas. Finally, for the modules where these tools don't give a conclusive answer, we use MAGMA algebra software (see \cite{BCP}) to compute the rank of the contravariant form $B$ on a certain graded piece of $M_c(\tau)$ and calculate the coefficients of the character formula from there.

We note that some of the results from this paper have been derived in other works by different methods. As mentioned above, \cite{E2} calculates which of the irreducible modules with trivial lowest weight representation are finite dimensional. We use this information in deriving the character formulas for them. \cite{R2}, Section 5.2.4, describes $L_{c}(\tau)$ in terms of $M_{c}(\tau)$ in case of ``blocks of defect one" - in our case $W=H_{3}$ these are $c=1/10$ and $c=1/6$. Let us also mention the paper \cite{M} that calculates the decomposition numbers for Iwahori-Hecke algebras associated to, among other groups, $H_{3}$. These numbers are related to the coefficients $n_{\tau,\sigma}$ in the Grothendieck group expression of $L_{c}(\tau)$ in terms of $M_{c}(\sigma)$. The results of \cite{M} can be derived easily from our results below. On the other hand, our results below don't follow directly, but could, with some work and in case $c\ne r/2$, be derived from the results of \cite{M}.

Potential further research would include calculating the characters for all the irreducible modules for the rational Cherednik algebra associated to $H_{4}$. This ought to be possible with similar methods, but with many more cases for case-by-case analysis, and additional programming difficulties coming from the size of $H_4$; or determining which of the irreducible representations described below of $H_{c}(H_3,\mathfrak{h})$ is unitary. 

The organization of the paper is as follows. Section 2 contains an overview of basic information about Cherednik algebras and their representations: definitions, basic properties, description of category $\mathcal{O}$ and standard and irreducible modules. It also contains facts about the group $H_3$ and its representations that we are going to use. Section 3 contains the statement of the main theorem. Section 4 presents a number of techniques we are going to use in the proof. Section 5 serves to reduce the number of parameters $c$ we consider, by quoting some previously known  equivalences of categories. As a consequence of this section, if we calculate the characters for $c\in \lbrace 1/10,1/6,1/5,1/3,1/2,3/2\rbrace$, we will know them for all $c\in \mathbb{C}$. In Sections 6-11 we do the main computational work of this paper, which is to describe the modules $L_{c}(\tau)$ in terms of $M_{c}(\sigma)$ for all $\tau$ and for  $c\in \lbrace 1/10,1/6,1/5,1/3,1/2,3/2\rbrace$.  

\subsection*{Acknowledgements} The authors are very grateful to Pavel
Etingof and for introducing us to this area of research, suggesting the problem, and devoting a lot of his time and energy to it through many helpful conversations. We thank Charles F. Dunkl for pointing out an error in the previous version of the text and to Stephen Griffeth for explaining to us Theorem \ref{equi1}. The work of both authors was supported by the Research Science Institute, and conducted in the Department of Mathematics at MIT. The work of M.B. was partially supported by the NSF grant DMS-0504847.

\section{ Preliminaries} 

\subsection{Rational Cherednik algebras}

Let $W$ be a finite group and $\mathfrak{h}$ a finite-dimensional faithful complex representation with a non-degenerate $W$-invariant inner product $\left< \cdot , \cdot\right> $. Denote by $S$ the set of reflections of $W$, meaning elements $s\in W$ such that $\mathrm{rk}(s-1)\mid_{\mathfrak{h}} = 1$. Since $W\subseteq O(\mathfrak{h})$, the only eigenvalue of $s$ different from $1$ is $-1$. $W$ also acts on $\mathfrak{h}^*$ and $s$ has the same eigenvalues there; denote by  $\alpha_s \in \mathfrak{h}^*$ and $\alpha_s^{\vee} \in \mathfrak{h}$ eigenvectors of $s$ with eigenvalues $-1$; chosen so that with respect to the natural pairing $(\cdot,\cdot): \mathfrak{h}^*\times \mathfrak{h}\to \mathbb{C}$ one has $( \alpha_s,\alpha_s^{\vee} ) = 2$. Let $c:S\to \mathbb{C}$ be an arbitrary conjugation invariant function. 

\begin{definition}
The rational Cherednik algebra $H_c(W,\mathfrak{h})$ is the quotient of $\mathbb{C}W \ltimes \mathrm{T}(\mathfrak{h}\oplus \mathfrak{h}^*)$ by the relations
$$ [x,x'] = 0,\quad  [y,y'] = 0 ,\quad [y,x] = ( y,x )  - \sum_{s \in S} c(s) ( \alpha_s,y )( x,\alpha_s^{\vee} ) s$$
for all $x, x' \in \mathfrak{h}^*$, $y, y' \in \mathfrak{h}$. 
\end{definition}

\begin{remark}Such an algebra can also be defined in the case there is no non-degenerate $W$-invariant inner product on $\mathfrak{h}$. In this case (complex) reflections $s$ can have eigenvalues $\lambda_{s}$ different from $-1$ and $1$. However, as this paper deals with a concrete Coxeter group, we shall use the above definition and keep assuming that the  form $\left< \cdot, \cdot \right> $ exists.
\end{remark}

\begin{remark} This algebra is sometimes denoted by $H_{1,c}(W,\mathfrak{h})$.\end{remark}

An analogue of the PBW theorem holds for $H_c(W,\mathfrak{h})$; meaning that as a vector space  $$H_c(W,\mathfrak{h})\cong \mathrm{S} \mathfrak{h}^* \otimes \mathbb{C}W \otimes \mathrm{S} \mathfrak{h}.$$

There is a very useful copy of $\mathfrak{sl}_2(\Bbb{C})$ in the algebra $H_c(W,\mathfrak{h})$. Fix an orthonormal basis $x_i$ in $\mathfrak{h}^*$  and the dual basis $y_i$ in $\mathfrak{h}$. Then this  copy of $\mathfrak{sl}_2(\Bbb{C})$ is spanned by $\mathbf{E} = \frac{1}{2}\sum_i x_i^2$, $\mathbf{F} =- \frac{1}{2}\sum_i y_i^2$, and a semisimple element 
$$\mathbf{h}=\sum_i x_iy_i+\frac{\dim \mathfrak{h}}{2}-\sum_{s\in S}c(s)s.$$ Direct calculation shows that $[\mathbf{h},x] = x, [\mathbf{h},y] = -y, [\mathbf{h},w] =0$ for $x \in \mathfrak{h}^*, y \in \mathfrak{h}, w\in W$, and from this it easily follows that $[\mathbf{h},\mathbf{E}] = 2\mathbf{E}, [\mathbf{h},\mathbf{F}] = -2\mathbf{F}, [\mathbf{E},\mathbf{F}] = \mathbf{h}$. It is a direct computation to check that this copy of $\mathfrak{sl}_{2}$ commutes with the elements of $H_{c}(W,\mathfrak{h})$ corresponding to the group $W$. The key to the usefulness of this subalgebra is that the element $\mathbf{h}$ acts as a grading element with respect to the grading given by $\mathrm{deg} x =1, \mathrm{deg} y =-1, \mathrm{deg} w=0, x \in \mathfrak{h}^*, y \in \mathfrak{h}, w\in W$. Thus, $H_c(W,\mathfrak{h})$ becomes an algebra with an inner grading and finite dimensional graded pieces which are representations of a finite group $W$. All the representations of $H_c(W,\mathfrak{h})$ that we will consider will also have this property.

\subsection{Category $\mathcal{O}$}

Define the category $\mathcal{O}_c(W,\mathfrak{h})$, or just $\mathcal{O}_c$, or just $\mathcal{O}$ (depending if we need to emphasize the algebra or if it is clear from the context) to be the category of $H_c(W,\mathfrak{h})$-modules which are finitely generated under the action of ${\rm S}\mathfrak{h}^*$ and locally nilpotent under the action of ${\rm S} \mathfrak{h}$. It's easy to see that $\mathbf{h}$ acts locally finitely on any module in this category. its generalized eigenspaces are finite dimensional representations of $W$. The category $\mathcal{O}_c$ contains all finite-dimensional $H_c(W,\mathfrak{h})$-modules. 

An important family of modules in category $\mathcal{O}_c$ are the standard or Verma modules $M_c(\tau)$. To define them, let 
$\tau$ be an irreducible representation of $W$. We call $\tau$ the lowest weight subspace of $M_c(\tau)$, and it plays the role of the one dimensional lowest weight vector space in Verma modules for Lie algebras. First define a structure of a module over $\Bbb{C}W \ltimes {\rm S}\mathfrak{h}\subseteq H_c(W,\mathfrak{h})$ on $\tau$ by letting ${\rm S}\mathfrak{h}$ act on $\tau$ by $0$. Then induce the action to the entire $H_c(W,\mathfrak{h})$, i.e. let ${\rm S}\mathfrak{h}^*$ act freely. So, as a vector space, the standard module $M_c(\tau)$ is isomorphic to ${\rm S}\mathfrak{h}^*\otimes \tau$. It is possible to make the action on this space a bit more explicit. It is clear how the elements of $\mathfrak{h}^*$ and $W$ act; and the action of elements of $y\in \mathfrak{h}$ can be described using the following Dunkl operators: $$D_{y}=\partial_{y}\otimes 1-\sum_{s\in S} c(s)\frac{(\alpha_s,y)}{\alpha_s}(1-s)\otimes s.$$

It is an easy direct computation to show that operators $\frac{1}{\alpha_s}(1-s)$ preserve ${\rm S}\mathfrak{h}$, and that the action of $D_y$ really coincides with the induced action of $y$.

Analogous to Lie theory, the sum of all the proper subrepresentations of $M_c(\tau)$, called $J_c(\tau)$, is the maximal proper submodule of $M_c(\tau)$. The quotient, called $L_c(\tau)$, is an irreducible module. All the irreducible modules in category $\mathcal{O}_{c}$ are of this form. The universal mapping properties analogous to the ones in Lie theory hold: for any module generated by the lowest weight representation $\tau$ (i.e., for any module $V$ that has a $W$ subrepresentation isomorphic to $\tau$, such that the action of ${\rm S}\mathfrak{h}$ on it is $0$ and the ${\rm S}\mathfrak{h}^*$ action on $\tau$ generates the entire $V$) there are unique surjective maps from $M_{c}(\tau)$ to $V$ and from $V$ to $L_{c}(\tau)$ that are identity on $\tau$.

As mentioned before, many things are unknown about the structure of $L_{c}(\tau)$. In this paper we will describe their structure for all the possible values of $c$ and $\tau$, for the exceptional Coxeter group $H_3$.

Whenever $c$ is clear from the context, such as in Sections \ref{ch1/10}-\ref{ch3/2}, we will write $L(\tau)$ for $L_{c}(\tau)$.

The grading element $\mathbf{h}=\sum_i x_iy_i+\frac{\dim \mathfrak{h}}{2}-\sum_{s\in S}c(s)s$ diagonalizes on $M_{c}(\tau)$. It acts on the lowest weight  $\tau$ by $h_{c}(\tau)=\frac{\dim \mathfrak{h}}{2}-\sum_{s\in S}c(s)s$. As $\sum_{s\in S}c(s)s$ is a central element of $\mathbb{C}W$, this is a constant depending on $c$ and $\tau$. If we then put a grading induced by $\mathbf{h}$ on $H_c(W,\mathfrak{h})\cong \mathrm{S} \mathfrak{h}^* \otimes \mathbb{C}W \otimes \mathrm{S} \mathfrak{h}$ i.e. $\mathrm{deg} x=1, x \in \mathfrak{h}^*, \mathrm{deg} y=-1, y \in \mathfrak{h}, \mathrm{deg} w=0, w\in W $, we see that $\mathbf{h}$ weights make $M_{c}(\tau)$ a graded representation. Denote the polynomials of degree $k$ by $\mathrm{S}^{k}\mathfrak{h}^*$ and the $j$ -th graded piece of $M_{c}(\tau)$ by $M_{c}(\tau)[j]$. The graded pieces of a standard module $M_{c}(\tau)$ are $\mathrm{S}^{k}\mathfrak{h}^*\otimes \tau \cong M_{c}(\tau)[k+h_{c}(\tau)]$. The dimension of this graded piece is ${k+\dim \mathfrak{h}^*-1 \choose \dim \mathfrak{h}^*-1} \cdot \dim \tau$, and it is $W$-invariant. This grading descends to all quotients of $M_{c}(\tau)$, most notably to $L_{c}(\tau)$. 

A simple but very useful observation is that the $\mathbf{h}$ weights are by definition the $\mathfrak{sl}_2$ weights that would appear in the decomposition of the module into $\mathfrak{sl}_2$ subrepresentations. So, if we are only interested in knowing which representations $V$ are finite dimensional, we can immediately put some obvious conditions on the weights that appear: all must be integral, the lowest one must be negative, and $\dim V[j]=\dim V[-j]$. 

One can define a contravariant symmetric bilinear form on $M_c(\tau)$, analogous to Shapovalov form on Verma modules in Lie theory. To start with, there is a $W$-invariant form $B$ on the lowest weight $\tau$. Extend it to $M_c(\tau)$ by requiring $B(x_i a,b)=B(a,y_ib), B(y_i a,b)=B(a,x_ib), B(wa,b)=B(a,wb)$ for all $a,b\in M_c(\tau), w\in W$, and $x_i, y_i$ orthonormal bases of $\mathfrak{h}^*, \mathfrak{h}$ dual to each other. It can be shown that this form is really well defined, that the different graded pieces of $M_c(\tau)$ are orthogonal to each other, and that the kernel of $B$ on $M_c(\tau)$ is exactly $J_c(\tau)$. This is very useful in computations; computing the dimension of kernel of $B$ in a certain graded piece of $M_{c}(\tau)$ can be easily done (by hand in lower and by computer in higher degrees), and gives us information about the size of $L_{c}(\tau)$. 

A module $V$ in category $\mathcal{O}$ is completely determined by its character. If $V=\oplus_{j}V[j]$, where $j\in \mathbf{C}$ is a generalized eigenvalue of $\mathbf{h}$ and $V[j] $ a generalized eigenspace, then  the character of $V$ is a function of $t\in \mathbb{C}$ and $w\in W$ given by $$\mathrm{ch}_{V}(w,t)= \mathrm {Tr}_{V}(wt^{\mathbf{h}}) =\sum_{j}t^{j} \mathrm {Tr}_{V[j] }(w).$$ It's easy to calculate the characters of the standard modules $M_{c}(\tau)$. Namely, if the character of a $W$ representation $\tau$ is $\chi_{\tau}$, and $\mathbf{h}$ acts on $\tau$ by a constant $h_{c}(\tau)$, then the character of $M_{c}(\tau)$ is given by $$\mathrm{ch}_{M_{c}(\tau)}(w,t)=\frac{\chi_{\tau}(w)t^{h_{c}(\tau)}}{\mathrm{det}_{\mathfrak{h}^*} (1-wt) }.$$

All modules in the category $\mathcal{O}$ have finite length. Since all the irreducible modules are $L_{c}(\tau)$ for all possible $\tau \in Irrep W$ (finite set), in the Grothendieck group $K_{0}(\mathcal{O}_{c})$ every module can be expressed as a finite linear combination with integer coefficients of modules $L_{c}(\tau)$. Express $M_{c}(\tau)$ like this, i.e. let $n'_{\tau,\sigma}$ be nonnegative integers such that in in Grothendieck group $$M_{c}(\tau)=\sum_{\sigma} n'_{\tau,\sigma} L_{c}(\sigma). \label{chMc}$$ A simple observation (see Lemma \ref{even}) is that $n'_{\tau,\tau}=1$, and, for $\tau \ne \sigma$, $n'_{\tau,\sigma}=0$ unless $h_c(\sigma)-h_c(\tau) $ is a positive integer. As a consequence of this, we can invert the matrix $[n'_{\tau,\sigma}]$ and it will still be upper triangular with integer (not necessarily positive) entries. In other words, we can find $n_{\tau,\sigma}\in \mathbb{Z}$, such that  $$L_{c}(\tau)=\sum_{\sigma} n_{\tau,\sigma} M_{c}(\sigma).$$ They will still satisfy $n_{\tau,\tau}=1$ and, for $\tau \ne \sigma$, $n_{\tau,\sigma}=0$ unless $h_c(\sigma)-h_c(\tau) $ is a positive integer. Of course, if such an expression holds for Grothendieck group, it will hold at the level of characters as well. As mentioned before, characters determine the irreducible modules completely, so if we calculate all the $n_{\tau,\sigma}$, and we know the formulas for characters of $M_{c}(\tau)$, we can consider the structure of $L_{c}(\tau)$ described. For example, then it is easy to determine which modules are finite dimensional. Finding all the possible $n_{\tau,\sigma}$ for all the possible values of $c$ and $\tau$ is what we do in this paper.

For references about Cherednik algebras and category $\mathcal{O}$, see \cite{E}, \cite{EM}, \cite{GGOR}.

\subsection{The group $H_3$}

We will study the case when $W=H_3$, the exceptional Coxeter group with the Coxeter graph 

\begin{center}
\begin{picture}(60,30)
\put(0,10){\line(1,0){60}}
\put(0,10){\circle*{5}}
\put(30,10){\circle*{5}}
\put(60,10){\circle*{5}}
\put(42,15){$5$}
\end{picture}
\end{center}

 It is the group of symmetries of the regular icosahedron, with Coxeter generators corresponding to the reflections along the planes with respective angles of $\pi/2$, $\pi/3$ and $\pi/5$ with each other. In this case, $\mathfrak{h}$ is the complexification of the $3$ dimensional real representation that realizes $H_3$ as such symmetry group. The scalar product that gives the structure of a Euclidean space to this $3$-dimensional real representation also gives an isomorphism $\mathfrak{h}\cong\mathfrak{h}^*$.

$H_3$ has $120$ elements and is isomorphic to $\mathbb{Z}_{2}\times A_{5}$. Here $\mathbb{Z}_2$ is a cyclic group of order $2$ containing the identity and the central symmetry of the icosahedron, and $A_5$ is a group of even permutations of the set of  $5$ elements, in this case the $5$ tetrahedra that are formed by centers of the faces of the icosahedron, which are permuted by rotations of the icosahedron. 

This presentation enables us to write the character table of $H_3$. The group $\mathbb{Z}_2$ has two one-dimensional irreducible representations, the trivial one and the signum one. The group $A_5$ has five irreducible representations: the trivial one, that we will call $\mathbf{1}$; a three dimensional one called $\mathbf{3}$, that realizes it as rotations of an icosahedron; another three dimensional one, called $\widetilde{\mathbf{3}}$, obtained from $\mathbf{3}$ by twisting by conjugation with the element $(12)\in S_5$ ($S_5$ is the symmetric group, and conjugating by $(12)$ in $S_5$ preserves $A_5\subset S_5$); a four dimensional representation $\mathbf{4}$ (the permutation representation of $A_5$ obtained from it acting on the $5$ tetrahedra is reducible; it has a $1$ dimensional trivial subrepresentation and $\mathbf{4}$ as irreducible components); and a $5$ dimensional $\mathbf{5}$, that is an irreducible subrepresentation of a $6$ dimensional representation arising from the fact that $A_5$ permutes the $6$ great diagonals of the icosahedron.  

Every irreducible representation of $H_3\cong \mathbb{Z}_2 \times A_5$ is a tensor product of an irreducible representation of $\mathbb{Z}_2$ and an irreducible representation of $A_5$. To simplify notation, let us denote the tensoring with the trivial representation by a subscript $+$, and tensoring with a signum representation by $-$. This makes sense because $+$ or $-$ now indicate whether the nontrivial element of $\mathbb{Z}_2$ acts by $1$ or by $-1$. So, we will denote $\mathrm{triv}\otimes \widetilde{\mathbf{3}}$ as $\widetilde{\mathbf{3}}_{+}$, or $\mathrm{sign}\otimes \mathbf{5}$ as $\mathbf{5}_{-}$. In the same style, write elements of $A_5$ as products of cycles, and elements of $\mathbb{Z}_2 \times A_5$ with a sign $+$ or $-$ in front to indicate which element of $\mathbb{Z}_2 $ is on their first coordinate (so, $(-1,(12)(34))=-(12)(34)$).

In this notation, $\mathfrak{h}\cong \mathfrak{h}^*\cong \mathbf{3}_{-}$.

The character table of $H_3$, in this notation, is Table \ref{H3CharTable}. For references about $H_3$ and its representations see \cite{H}, \cite{FH}.

\begin{table}[h!]
\begin{center}
\begin{tabular}{|c||c|c||c|c||c|c||c|c||c|c|} \hline
& Id & -Id & (123) & -(123) & (12)(34) & -(12)(34) & (12345) & -(12345) & (13245) & -(13245)\\ \hline
\# & 1 & 1 & 20 & 20 & 15 & 15 & 12 & 12 & 12 & 12 \\ \hline
$\mathbf{1}_{+}$ & 1 & 1 & 1 & 1 & 1 & 1 & 1 & 1 & 1 & 1\\ \hline
$\mathbf{1}_{-}$ & 1 & -1 & 1 & -1 & 1 & -1 & 1 & -1 & 1 & -1\\ \hline
$\mathbf{3}_{+}$ & 3 & 3 & 0 & 0 & -1 & -1 & $\frac{1+\sqrt{5}}{2}$ & $\frac{1+\sqrt{5}}{2}$ & $\frac{1-\sqrt{5}}{2}$ & $\frac{1-\sqrt{5}}{2}$\\ \hline
$\mathbf{3}_{-}$ & 3 & -3 & 0 & 0 & -1 & 1 & $\frac{1+\sqrt{5}}{2}$ & $\frac{-1-\sqrt{5}}{2}$ & $\frac{1-\sqrt{5}}{2}$ & $\frac{-1+\sqrt{5}}{2}$\\ \hline
$\widetilde{\mathbf{3}}_+$ & 3 & 3 & 0 & 0 & -1 & -1 & $\frac{1-\sqrt{5}}{2}$ & $\frac{1-\sqrt{5}}{2}$ & $\frac{1+\sqrt{5}}{2}$ & $\frac{1+\sqrt{5}}{2}$\\ \hline
$\widetilde{\mathbf{3}}_-$ & 3 & -3 & 0 & 0 & -1 & 1 & $\frac{1-\sqrt{5}}{2}$ & $\frac{-1+\sqrt{5}}{2}$ & $\frac{1+\sqrt{5}}{2}$ & $\frac{-1-\sqrt{5}}{2}$\\ \hline
$\mathbf{4}_{+}$ & 4 & 4 & 1 & 1 & 0 & 0 & -1 & -1 & -1 & -1\\ \hline
$\mathbf{4}_{-}$ & 4 & -4 & 1 & -1 & 0 & 0 & -1 & 1 & -1 & 1\\ \hline
$\mathbf{5}_{+}$ & 5 & 5 & -1 & -1 & 1 & 1 & 0 & 0 & 0 & 0\\ \hline
$\mathbf{5}_{-}$ & 5 & -5 & -1 & 1 & 1 & -1 & 0 & 0 & 0 & 0\\ \hline
\end{tabular}
\end{center}
\caption{The character table for $H_3\cong \Bbb{Z}_2 \times A_5$}
\label{H3CharTable}
\end{table}

There is only one conjugacy class of reflections in $H_3$, namely the class of $-(12)(34)$, with $15$ reflections in it. Since $c:S\to \mathbb{C}$ needs to be a conjugation invariant function, in the case of $H_3$ it is just a complex constant. 

We can also easily calculate the action of the central element $\sum_{s\in S} s$ on any representation. For example, in $\mathbf{5}_{-}$, it is a constant on a $5$ dimensional space, whose trace is $\mathrm{tr}=\sum_{s\in S} \mathrm{tr} s=-15$, so it is $-15/5=3$. Doing this calculation for every irreducible representation $\tau$, we get Table \ref{const sum s}.

\begin{table}[h!]
\begin{center}
\begin{tabular}{|c|c|c|c|c|c|c|c|c|c|} \hline
$\mathbf{1}_{+}$ & $\mathbf{1}_{-}$ & $\mathbf{3}_{+}$ & $\mathbf{3}_{-}$ & $\widetilde{\mathbf{3}}_+$ & $\widetilde{\mathbf{3}}_-$ & $\mathbf{4}_{+}$ & $\mathbf{4}_{-}$ & $\mathbf{5}_{+}$ & $\mathbf{5}_{-}$ \\ \hline
15 & -15 & -5 & 5 & -5 & 5 & 0 & 0 & 3 & -3 \\ \hline
\end{tabular}
\end{center}
\caption{The action of the central element $\sum_{s\in S} s\in H_3$ on all $\tau$}
\label{const sum s}
\end{table}

Table  \ref{const sum s} now enables us to calculate the action of $\mathbf{h}$ on any lowest weight $\tau$, as $h_{c}(\tau)=\frac{3}{2}-c\sum_{s\in S}s|_\tau$.

\section{Main theorem}

\begin{theorem}
For the Coxeter group $H_3$, its reflection representation $\mathfrak{h}$, $c$ any complex number, and $\tau$ an irreducible representation of $H_3$, the expression in the Grothendieck group  $K_{0}(\mathcal{O}_{c})$ for the irreducible module $L_{c}(\tau)$ in terms of standard modules $M_{c}(\tau)$ is as below.  Any module $L_{c}(\tau)$ for which we do not explicitly write its dimension is infinite dimensional. We leave out the index $c$ in $L_{c}(\tau)$ and $M_{c}(\tau)$ whenever it is clear from the context. Here $r\in \mathbb{N}, d\in \{ 2,3,5,6,10\}$, and all fractions $r/d$ are reduced.

\begin{itemize}
\item If $c$ is not of the form $c=r/d$ or $c=-r/d$, then for all $\tau$, $$L_{c}(\tau)=M_{c}(\tau).$$

\item If $c=-r/d$, then the formulas for $L_{c}(\tau)$ in terms of $M_{c}(\sigma)$ follow from formulas for $L_{-c}(\mathbf{1}_{-}\otimes \tau)$ in terms of $M_{c}(\mathbf{1}_{-}\otimes \sigma)$, which are given below. More precisely, if $$L_{r/d}(\tau)=\sum_{\sigma}n_{\tau,\sigma}M_{r/d}(\sigma)$$ then  $$L_{-r/d}(\mathbf{1}_{-}\otimes \tau)=\sum_{\sigma}n_{\tau,\sigma}M_{-r/d}(\mathbf{1}_{-}\otimes \sigma).$$ Consequently,  $L_{c}(\tau)$ is finite dimensional if and only if $L_{-c}(\mathbf{1}_{-}\otimes \tau)$ is.

\item $\mathbf{c=r/10, r\ne 3,7 \, (\mathrm{mod}\, 10)}$
\begin{eqnarray*}
L(\mathbf{1}_{+}) & = & M(\mathbf{1}_{+})- M(\mathbf{3}_{-})+M(\mathbf{3}_{+})-M(\mathbf{1}_{-}) \\
L(\mathbf{3}_{+})& = &M(\mathbf{3}_{+})-M(\mathbf{1}_{-})\\
L(\mathbf{3}_{-})& = &M(\mathbf{3}_{-})-M(\mathbf{3}_{+})+M(\mathbf{1}_{-})
\end{eqnarray*}

Every $L_{r/10}(\mathbf{1}_{+})$ is finite dimensional, with $\dim L_{r/10}(\mathbf{1}_{+})=r^3$ and  $$\mathrm{ch}_{L_{r/10}(\mathbf{1}_{+})}(w,t)=\frac{\det_{\mathfrak{h}^*}(1-wt^r)}{\det_{\mathfrak{h}^*}(1-wt)}.$$

\item $\mathbf{c=r/10, r= 3,7 \, (\mathrm{mod}\, 10)}$
\begin{eqnarray*}
L(\mathbf{1}_{+}) & = & M(\mathbf{1}_{+})- M(\tilde{\mathbf{3}}_{-})+M(\tilde{\mathbf{3}}_{+})-M(\mathbf{1}_{-}) \\
L(\tilde{\mathbf{3}}_{+})& = &M(\tilde{\mathbf{3}}_{+})-M(\mathbf{1}_{-})\\
L(\tilde{\mathbf{3}}_{-})& = &M(\tilde{\mathbf{3}}_{-})-M(\tilde{\mathbf{3}}_{+})+M(\mathbf{1}_{-})
\end{eqnarray*}

Every $L_{r/10}(\mathbf{1}_{+})$ is finite dimensional, with $\dim L_{r/10}(\mathbf{1}_{+})=r^3$.

\item $\mathbf{c=r/6}$
\begin{eqnarray*}
L(\mathbf{1}_{+})&=&M(\mathbf{1}_{+})-M(\mathbf{5}_{+})+M(\mathbf{5}_{-})-M(\mathbf{1}_{-})\\
L(\mathbf{5}_{+})&=&M(\mathbf{5}_{+})-M(\mathbf{5}_{-})+M(\mathbf{1}_{-})\\
L(\mathbf{5}_{-})&=&M(\mathbf{5}_{-})-M(\mathbf{1}_{-})
\end{eqnarray*}
Every $L_{r/6}(\mathbf{1}_{+})$is finite dimensional, with $\dim L_{r/6}(\mathbf{1}_{+})=5r^3$ and  $$\mathrm{ch}_{L_{r/6}(\mathbf{1}_{+})}=\frac{\det_{\mathfrak{h}^*}(1-wt^r)}{\det_{\mathfrak{h}^*}(1-wt)}\cdot \big( \chi_{\mathbf{1}_{+}}t^{-r} + \chi_{\mathbf{3}_{-}}+\chi_{\mathbf{1}_{+}}t^r \big).$$

\item $\mathbf{c=r/5, r= 1,9 \, (\mathrm{mod}\, 10)}$
\begin{eqnarray*}
L( \mathbf{1}_{+})&=&M( \mathbf{1}_{+})-M( \mathbf{4}_{-})+M( \widetilde{\mathbf{3}}_{+})\\
L(\widetilde{\mathbf{3}}_{-})&=&M(\widetilde{\mathbf{3}}_{-})-M(\mathbf{4}_{+})+M(\mathbf{1}_{-})\\
L(\mathbf{4}_{+})&=&M(\mathbf{4}_{+})-M(\mathbf{1}_{-})\\
L(\mathbf{4}_{-})&=&M(\mathbf{4}_{-})-M(\widetilde{\mathbf{3}}_{+})\\
\end{eqnarray*}

\item $\mathbf{c=r/5, r= 2,8 \, (\mathrm{mod}\, 10)}$
\begin{eqnarray*}
L( \mathbf{1}_{+})&=&M( \mathbf{1}_{+})-M( \mathbf{4}_{+})+M(\mathbf{3}_{+})\\
L(\mathbf{3}_{-})&=&M(\mathbf{3}_{-})-M(\mathbf{4}_{-})+M(\mathbf{1}_{-})\\
L(\mathbf{4}_{-})&=&M(\mathbf{4}_{-})-M(\mathbf{1}_{-})\\
L(\mathbf{4}_{+})&=&M(\mathbf{4}_{+})-M(\mathbf{3}_{+})\\
\end{eqnarray*}

\item $\mathbf{c=r/5, r= 3,7 \, (\mathrm{mod}\, 10)}$
\begin{eqnarray*}
L( \mathbf{1}_{+})&=&M( \mathbf{1}_{+})-M( \mathbf{4}_{-})+M(\mathbf{3}_{+})\\
L(\mathbf{3}_{-})&=&M(\mathbf{3}_{-})-M(\mathbf{4}_{+})+M(\mathbf{1}_{-})\\
L(\mathbf{4}_{+})&=&M(\mathbf{4}_{+})-M(\mathbf{1}_{-})\\
L(\mathbf{4}_{-})&=&M(\mathbf{4}_{-})-M(\mathbf{3}_{+})\\
\end{eqnarray*}

\item $\mathbf{c=r/5, r= 4,6 \, (\mathrm{mod}\, 10)}$
\begin{eqnarray*}
L( \mathbf{1}_{+})&=&M( \mathbf{1}_{+})-M( \mathbf{4}_{+})+M( \widetilde{\mathbf{3}}_{+})\\
L(\widetilde{\mathbf{3}}_{-})&=&M(\widetilde{\mathbf{3}}_{-})-M(\mathbf{4}_{-})+M(\mathbf{1}_{-})\\
L(\mathbf{4}_{-})&=&M(\mathbf{4}_{-})-M(\mathbf{1}_{-})\\
L(\mathbf{4}_{+})&=&M(\mathbf{4}_{+})-M(\widetilde{\mathbf{3}}_{+})\\
\end{eqnarray*}

\item $\mathbf{c=r/3, r}$ \textbf{odd}
\begin{eqnarray*}
L(\mathbf{1}_{+})&=&M(\mathbf{1}_{+})-M(\mathbf{5}_{+})+M(\mathbf{4}_{-})\\
L(\mathbf{4}_{+})&=&M(\mathbf{4}_{+})-M(\mathbf{5}_{-})+M(\mathbf{1}_{-})\\
L(\mathbf{5}_{-})&=& M(\mathbf{5}_{-})-M(\mathbf{1}_{-})\\
L(\mathbf{5}_{+})&=&M(\mathbf{5}_{+})-M(\mathbf{4}_{-})
\end{eqnarray*}

\item $\mathbf{c=r/3, r}$ \textbf{even}
\begin{eqnarray*}
L(\mathbf{1}_{+})&=&M(\mathbf{1}_{+})-M(\mathbf{5}_{+})+M(\mathbf{4}_{+})\\
L(\mathbf{4}_{-})&=&M(\mathbf{4}_{-})-M(\mathbf{5}_{-})+M(\mathbf{1}_{-})\\
L(\mathbf{5}_{-})&=& M(\mathbf{5}_{-})-M(\mathbf{1}_{-})\\
L(\mathbf{5}_{+})&=&M(\mathbf{5}_{+})-M(\mathbf{4}_{+})
\end{eqnarray*}

\item $\mathbf{c=r/2}$
\begin{eqnarray*}
L(\mathbf{1}_{+})&=& M(\mathbf{1}_{+})-M(\mathbf{3}_{-})-M(\widetilde{\mathbf{3}}_{-})+M(\mathbf{5}_{+})-M(\mathbf{5}_{-})+M(\mathbf{3}_{+})+M(\widetilde{\mathbf{3}}_{+})-M(\mathbf{1}_{-})\\
L(\mathbf{3}_{+})&=&M(\mathbf{3}_{+})-M(\mathbf{1}_{-})\\
L(\mathbf{3}_{-})&=& M(\mathbf{3}_{-})-M(\mathbf{5}_{+})+M(\mathbf{5}_{-})-M(\mathbf{3}_{+})\\
L(\widetilde{\mathbf{3}}_{+})&=&M(\widetilde{\mathbf{3}}_{+})-M(\mathbf{1}_{-})\\
L(\widetilde{\mathbf{3}}_{-})&=& M(\widetilde{\mathbf{3}}_{-})-M(\mathbf{5}_{+})+M(\mathbf{5}_{-})-M(\widetilde{\mathbf{3}}_{+})\\
L(\mathbf{5}_{+})&=& M(\mathbf{5}_{+})-2\cdot M(\mathbf{5}_{-})+M(\mathbf{3}_{+})+M(\widetilde{\mathbf{3}}_{+})-M(\mathbf{1}_{-}).\\
L(\mathbf{5}_{-})&=&M(\mathbf{5}_{-})-M(\mathbf{3}_{+})-M(\widetilde{\mathbf{3}}_{+})+M(\mathbf{1}_{-})
\end{eqnarray*}
For every $r$, three of these modules are finite dimensional, with $\dim L_{r/2}(\mathbf{1}_{+})=115r^3$, $\dim L_{r/2}(\mathbf{3}_{-})=10r^3$, and $\dim L_{r/2}(\widetilde{\mathbf{3}}_{-})=10r^3$, and $$ch_{L_{r/2}(\mathbf{3}_{-})}(w,t)=\frac{\det_{\mathfrak{h}^*}(1-wt^r)}{\det_{\mathfrak{h}^*}(1-wt)}\cdot \big( \chi_{\mathbf{3}_{-}}t^{-r}+\chi_{\mathbf{1}_{+}}+\chi_{\mathbf{3}_{+}}+\chi_{\mathbf{3}_{-}}t^r \big),$$ $$ch_{L_{r/2}(\widetilde{\mathbf{3}}_{-})}(w,t)=\frac{\det_{\mathfrak{h}^*}(1-wt^r)}{\det_{\mathfrak{h}^*}(1-wt)}\cdot  \big( \chi_{\widetilde{\mathbf{3}}_{-}}t^{-r}+\chi_{\mathbf{4}_{+}}+\chi_{\widetilde{\mathbf{3}}_{-}}t^r \big) .$$

\end{itemize}
\label{main}
\end{theorem}

\begin{proof}
To see that the only values of $c$ for which the above formulas are nontrivial are $c=r/d$ with $r,d$ as above, see Section \ref{rationalc}. For a proof that only $c>0$ need to be considered, see Sections \ref{c ne 0} and \ref{c>0}. Theorems \ref{1/10}, \ref{1/6}, \ref{1/5}, \ref{1/3} calculate the Grothendieck group expressions for $c=1/d, d\in \{ 10,6,5,3\}$, and Theorem \ref{equi1} lets us conclude the above formulas for all $c=r/d, d\in \{ 10,6,5,3\}, r>0$. Similarly, Theorems \ref{1/2} and \ref{3/2} calculate the Grothendieck group expressions for $c=1/2$ and $c=3/2$, while Theorem \ref{equi2} and Lemma \ref{equi3} allow us to derive the above formulas for $c=r/d, r>3$.
\end{proof}

For completeness and reader's convenience, we also calculate the Grothendieck group expressions of $M_{c}(\tau)$ in terms of $L_{c}(\tau)$; i.e. calculate which irreducible modules appear in the composition series of standard modules. 

\begin{theorem}
For the Coxeter group $H_3$, its reflection representation $\mathfrak{h}$, $c$ any complex number, and $\tau$ an irreducible representation of $H_3$, the expression in the Grothendieck group  $K_{0}(\mathcal{O}_{c})$ for the standard module $M_{c}(\tau)$ in terms of irreducible modules $L_{c}(\tau)$ is as below.  We leave out the index $c$ in $L_{c}(\tau)$ and $M_{c}(\tau)$ whenever it is clear from the context. Here $r\in \mathbb{N}, d\in \{ 2,3,5,6,10\}$, and all fractions $r/d$ are reduced.

\begin{itemize}

\item If $c$ is not of the form $c=r/d$ or $c=-r/d$, then for all $\tau$, $$M_{c}(\tau)=L_{c}(\tau).$$

\item If $c=-r/d$, then the formulas for $M_{c}(\tau)$ in terms of $L_{c}(\sigma)$ follow from formulas for $M_{-c}(\mathbf{1}_{-}\otimes \tau)$ in terms of $L_{c}(\mathbf{1}_{-}\otimes \sigma)$, which are given below. More precisely, if $$M_{r/d}(\tau)=\sum_{\sigma}n'_{\tau,\sigma}L_{r/d}(\sigma)$$ then  $$M_{-r/d}(\mathbf{1}_{-}\otimes \tau)=\sum_{\sigma}n'_{\tau,\sigma}L_{-r/d}(\mathbf{1}_{-}\otimes \sigma).$$

\item $\mathbf{c=r/10, r\ne 3,7 \, (\mathrm{mod}\, 10)}$
\begin{eqnarray*}
M(\mathbf{1}_{+}) & = & L(\mathbf{1}_{+})+ L(\mathbf{3}_{-})\\
M(\mathbf{3}_{+})& = &L(\mathbf{3}_{+})+L(\mathbf{1}_{-})\\
M(\mathbf{3}_{-})& = &L(\mathbf{3}_{-})+L(\mathbf{3}_{+})
\end{eqnarray*}

\item $\mathbf{c=r/10, r= 3,7 \, (\mathrm{mod}\, 10)}$
\begin{eqnarray*}
M(\mathbf{1}_{+}) & = & L(\mathbf{1}_{+})+ L(\tilde{\mathbf{3}}_{-})\\
M(\tilde{\mathbf{3}}_{+})& = &L(\tilde{\mathbf{3}}_{+})+L(\mathbf{1}_{-})\\
M(\tilde{\mathbf{3}}_{-})& = &L(\tilde{\mathbf{3}}_{-})+L(\tilde{\mathbf{3}}_{+})
\end{eqnarray*}

\item $\mathbf{c=r/6}$
\begin{eqnarray*}
M(\mathbf{1}_{+})&=&L(\mathbf{1}_{+})+L(\mathbf{5}_{+})\\
M(\mathbf{5}_{+})&=&L(\mathbf{5}_{+})+L(\mathbf{5}_{-})\\
M(\mathbf{5}_{-})&=&L(\mathbf{5}_{-})+L(\mathbf{1}_{-})
\end{eqnarray*}

\item $\mathbf{c=r/5, r= 1,9 \, (\mathrm{mod}\, 10)}$
\begin{eqnarray*}
M( \mathbf{1}_{+})&=&L( \mathbf{1}_{+})+L( \mathbf{4}_{-})\\
M(\widetilde{\mathbf{3}}_{-})&=&L(\widetilde{\mathbf{3}}_{-})+L(\mathbf{4}_{+})\\
M(\mathbf{4}_{+})&=&L(\mathbf{4}_{+})+L(\mathbf{1}_{-})\\
M(\mathbf{4}_{-})&=&L(\mathbf{4}_{-})+L(\widetilde{\mathbf{3}}_{+})
\end{eqnarray*}

\item $\mathbf{c=r/5, r= 2,8 \, (\mathrm{mod}\, 10)}$
\begin{eqnarray*}
M( \mathbf{1}_{+})&=&L( \mathbf{1}_{+})+L( \mathbf{4}_{+})\\
M(\mathbf{3}_{-})&=&L(\mathbf{3}_{-})+L(\mathbf{4}_{-})\\
M(\mathbf{4}_{-})&=&L(\mathbf{4}_{-})+L(\mathbf{1}_{-})\\
M(\mathbf{4}_{+})&=&L(\mathbf{4}_{+})+L(\mathbf{3}_{+})
\end{eqnarray*}

\item $\mathbf{c=r/5, r= 3,7 \, (\mathrm{mod}\, 10)}$
\begin{eqnarray*}
M( \mathbf{1}_{+})&=&L( \mathbf{1}_{+})+L( \mathbf{4}_{-})\\
M(\mathbf{3}_{-})&=&L(\mathbf{3}_{-})+L(\mathbf{4}_{+})\\
M(\mathbf{4}_{+})&=&L(\mathbf{4}_{+})+L(\mathbf{1}_{-})\\
M(\mathbf{4}_{-})&=&L(\mathbf{4}_{-})+L(\mathbf{3}_{+})
\end{eqnarray*}

\item $\mathbf{c=r/5, r= 4,6 \, (\mathrm{mod}\, 10)}$
\begin{eqnarray*}
M( \mathbf{1}_{+})&=&L( \mathbf{1}_{+})+L( \mathbf{4}_{+})\\
M(\widetilde{\mathbf{3}}_{-})&=&L(\widetilde{\mathbf{3}}_{-})+L(\mathbf{4}_{-})\\
M(\mathbf{4}_{-})&=&L(\mathbf{4}_{-})+L(\mathbf{1}_{-})\\
M(\mathbf{4}_{+})&=&L(\mathbf{4}_{+})+L(\widetilde{\mathbf{3}}_{+})
\end{eqnarray*}

\item $\mathbf{c=r/3, r}$ \textbf{odd}
\begin{eqnarray*}
M(\mathbf{1}_{+})&=&L(\mathbf{1}_{+})+L(\mathbf{5}_{+})\\
M(\mathbf{4}_{+})&=&L(\mathbf{4}_{+})+L(\mathbf{5}_{-})\\
M(\mathbf{5}_{-})&=& L(\mathbf{5}_{-})+L(\mathbf{1}_{-})\\
M(\mathbf{5}_{+})&=&L(\mathbf{5}_{+})+L(\mathbf{4}_{-})
\end{eqnarray*}

\item $\mathbf{c=r/3, r}$ \textbf{even}
\begin{eqnarray*}
M(\mathbf{1}_{+})&=&L(\mathbf{1}_{+})+L(\mathbf{5}_{+})\\
M(\mathbf{4}_{-})&=&L(\mathbf{4}_{-})+L(\mathbf{5}_{-})\\
M(\mathbf{5}_{-})&=& L(\mathbf{5}_{-})+L(\mathbf{1}_{-})\\
M(\mathbf{5}_{+})&=&L(\mathbf{5}_{+})+L(\mathbf{4}_{+})
\end{eqnarray*}

\item $\mathbf{c=r/2}$
\begin{eqnarray*}
M(\mathbf{1}_{+})&=& L(\mathbf{1}_{+})+L(\mathbf{3}_{-})+L(\widetilde{\mathbf{3}}_{-})+L(\mathbf{5}_{+})+L(\mathbf{5}_{-})+L(\mathbf{1}_{-})\\
M(\mathbf{3}_{+})&=&L(\mathbf{3}_{+})+L(\mathbf{1}_{-})\\
M(\mathbf{3}_{-})&=& L(\mathbf{3}_{-})+L(\mathbf{5}_{+})+L(\mathbf{5}_{-})+L(\mathbf{3}_{+})+L(\mathbf{1}_{-})\\
M(\widetilde{\mathbf{3}}_{+})&=&L(\widetilde{\mathbf{3}}_{+})+L(\mathbf{1}_{-})\\
M(\widetilde{\mathbf{3}}_{-})&=& L(\widetilde{\mathbf{3}}_{-})+L(\mathbf{5}_{+})+L(\mathbf{5}_{-})+L(\widetilde{\mathbf{3}}_{+})+L(\mathbf{1}_{-})\\
M(\mathbf{5}_{+})&=& L(\mathbf{5}_{+})+2\cdot L(\mathbf{5}_{-})+L(\mathbf{3}_{+})+L(\widetilde{\mathbf{3}}_{+})+L(\mathbf{1}_{-}).\\
M(\mathbf{5}_{-})&=&L(\mathbf{5}_{-})+L(\mathbf{3}_{+})+L(\widetilde{\mathbf{3}}_{+})+L(\mathbf{1}_{-})
\end{eqnarray*}

\end{itemize} 
\end{theorem}

\section{Several techniques}
In this section we state several previously known results that we use in our computations.

We will use the following lemma several times (see \cite{ES}, lemma 3.5):  
\begin{lemma} Let $\sigma\subseteq \mathfrak{h}^*\otimes \tau=\mathrm{S}^1 \mathfrak{h}^*\otimes \tau \subseteq \mathrm{S} \mathfrak{h}^*\otimes \tau \cong M_c(\tau)$ be an irreducible $W$-subrepresentation. Then the elements of $\mathfrak{h}$ act on $\sigma$ by $0$
(i.e. $\sigma$ consists of singular vectors) if and only if $$h_c(\sigma)-h_c(\tau)=1.$$
\label{Lemma35}
\end{lemma}

The following lemma is a slightly stronger version of Corollary 2.20 in \cite{GGOR}, valid in case of $H_3$.
\begin{lemma}
 $n_{\tau,\sigma}=0$ unless $h_{c}(\sigma)-h_{c}(\tau)$ is a positive integer. If $-\mathrm{Id}\in H_3$ acts of both $\tau$ and $\sigma$ by the same constant (either $1$ or $-1$), this integer must be even; otherwise this integer must be odd.
\label{even}
\end{lemma}

\begin{proof}
 As the matrix of integers $[n_{\tau,\sigma}]$ is an inverse of the matrix $[n'_{\tau,\sigma}]$, it is enough to prove this for positive integers $n'_{\tau,\sigma}$ (because the matrix $[n'_{\tau,\sigma}]$, with an appropriate ordering on $\tau$, becomes not only upper triangular with $1$ on diagonal, but also block diagonal with $\tau$ and $\sigma$ being in the same block if $h_{c}(\sigma)-h_{c}(\tau)$ is an integer). They are the coefficients in the composition series of $M_{c}(\tau)=\sum_{\sigma}n'_{\tau,\sigma}L_{c}(\sigma)$. If $L_c(\sigma)$ appears in the composition series of $M_{c}(\tau)$, then $\sigma$ is a $W$ subrepresentation of some degree $k>0$ graded piece $\mathrm{S}^k\mathfrak{h}^{*}\otimes \tau$. So, the action of $\mathbf{h}$ on $\sigma$ must at the same time be $h_{c}(\sigma)$ and $h_{c}(\tau)+k$, which means $h_{c}(\sigma)-h_{c}(\tau)=k$. 
 
 Now consider the action of an element $-\mathrm{Id}\in W$ on graded pieces. It is a central element acting by $+1$ or $-1$ on any irreducible representation. It acts by $-1$ on $\mathfrak{h}^*=\mathbf{3}_{-}$; so it acts by $(-1)^k$ on $\mathrm{S}^k\mathfrak{h}^{*}$. If $L_{c}(\sigma)$ appears in the composition series of $M_{c}(\tau)$, then it must act on both $\sigma$ and $\mathrm{S}^k\mathfrak{h}^{*}\otimes \tau$ the same. This proves the second claim of the lemma.
 \end{proof}

\subsection{Support of a module}
We will also use the main result from the recent paper \cite{E2}. Let $c=r/d\in\mathbb{Q}_{+}$ be a rational parameter. This paper considers modules $L_{c}(\mathbb{C_+})$. Such a module is a quotient of $M_{c}(\mathbb{C_+})\cong \mathrm{S}\mathfrak{h}^*\cong \mathbb{C}[\mathfrak{h}]$, so one can consider their structure as a module over the ring $\mathbb{C}[\mathfrak{h}]$. In this language $J_{c}(\mathbb{C_+})$ is an ideal $L_{c}(\mathbb{C_+})$, and $L_{c}(\mathbb{C_+})=\mathbb{C}\mathfrak{h}/J_{c}(\mathbb{C_+})$. One can then consider its support in the sense of commutative algebra. For $a\in \mathfrak{h}$ denote by $W_a$ the stabilizer of $a$ in $W$. Let $d_i=d_i(W)$ be the degrees of basic invariants, i.e. degrees of the generators of the polynomial algebra $\mathbb{C}[\mathfrak{h}]^W$. Let $l(w)$ be the length of $w\in W$. Define the Poincar\'{e} polynomial of $W$ to be $$P_{W}(q)=\sum_{w\in W}q^{l(w)}=\prod_{i}\frac{1-q^{d_i}}{1-q}.$$
Theorem 3.1. in \cite{E2} then states:
\begin{theorem}
 A point $a\in\mathfrak{h}$ belongs to the support of $L_c$ if and only if $$\frac{P_W}{P_{W_a}}(e^{2\pi i c})\ne 0$$
 i.e. if and only if 
$$\# \lbrace i | d \text{ divides } d_i(W)  \rbrace= \# \lbrace i | d \text{ divides } d_i(W_a)  \rbrace.$$ 
\label{support}
\end{theorem}

Geometry of the support of $L_{c}(\mathbf{1}_{+})$ tells us a lot. We can regard $L_{c}(\mathbf{1}_{+})=\mathbb{C}[\mathfrak{h}]/J_{c}(\mathbf{1}_{+})$ as a ring, and then its support is just the subvariety of $\mathfrak{h}$ determined by the ideal $J_{c}(\mathbf{1}_{+})$, i.e. $\mathrm{Spec} L_{c}(\mathbf{1}_{+})$. So, its dimension is equal to the degree of the pole at $t=1$ of the Hilbert-Poincar\' e series of $L_{c}(\mathbf{1}_{+})$ with respect to the usual grading on $\mathbb{C}[\mathfrak{h}]$. This grading and the grading by $\mathbf{h}$ action differ by a constant $h_c(\mathbf{1}_{+})$, so Hilbert-Poincar\' e series defined using these two gradings differ by a constant power of $t$, so they have the same order of pole at $t=1$. But Hilbert-Poincar\' e series of $L_{c}(\mathbf{1}_{+})$ with respect to the $\mathbf{h}$ grading is just the character of $L_{c}(\mathbf{1}_{+})$ evaluated at $w=\mathrm{Id}$. This will help us determine the coefficients in the character formulas for $L_{c}(\mathbf{1}_{+})$. 

In particular, $L_{c}(\mathbb{C}_{+})$ is finite dimensional if and only if there is no pole, meaning if the support is zero dimensional (in this case, it has to be $0$, but that will not be important to us).

\subsection{Parabolic induction and restriction functors}\label{indchap}
We use the results of \cite{BE}. For $W'$ a parabolic subgroup of $W$, for $b\in \mathfrak{h}$ such that $W'$ is the stabilizer of $b$ in $W$, and for $\mathfrak{h}'=\mathfrak{h}/\mathfrak{h}^{W'}$ reflection representation of $W'$, this paper defines functors $\Res_{b}$ and $\Ind_{b}$, depending on the point $b$. For simplicity assume $c\in\mathbb{C}$ is a constant, as it is in our case. Then these functors are between categories $\mathcal{O}$ for $H_{c}(W,\mathfrak{h})$ and for $H_{c}(W',\mathfrak{h}')$:  $$\Res_{b}: \mathcal{O}_{c}(W,\mathfrak{h})\to \mathcal{O}_{c}(W',\mathfrak{h}')$$ $$\Ind_{b}: \mathcal{O}_{c}(W',\mathfrak{h}')\to \mathcal{O}_{c}(W,\mathfrak{h}).$$ We will not use their construction nor their deeper properties, but just the following results (Proposition 3.14 from \cite{BE}):

\begin{proposition}\label{ind}
At the level of Grothendieck groups, $\Res_{b}$ and $\Ind_{b}$ applied to standard modules $L_{c}(\tau)$ behave like the induction and restriction functors for finite groups $W$ and $W'$; namely, for $\tau \in \hat{W}, \sigma \in \hat{W}'$:
$$\Res_{b}(M_{c}(\tau))=\sum_{\xi \in \hat{W'}} (\dim \Hom_{W'}(\xi,\tau))\cdot M_{c}(\xi)$$
$$\Ind_{b}(M_{c}(\sigma))=\sum_{\xi \in \hat{W}} (\dim \Hom_{W'}(\sigma,\xi))\cdot M_{c}(\xi)$$
\end{proposition}

\section{ Equivalences of categories}

In this section, we gather results from various papers to reduce the number of parameters $c$ for which we need to do the calculations. The idea is that there are many equivalences of categories between $\mathcal{O}_c(W,\mathfrak{h})$ for various $c$, for which we know where they map standard and irreducible modules. So calculating of characters of $L_{c}(\tau)$ for a specific $c$ and then applying the appropriate functor gives us character formulas for the whole family of parameters $c$. Also, for many $c$ it is known that the category is semisimple, so there is nothing to compute there. 

\subsection{WLOG $c\ne 0$}\label{c ne 0}

In case $c=0$, the grading element $\mathbf{h}$ acts on all the $\tau$ by the same constant, $h_{0}(\tau)=3/2$. Since $n_{\tau,\sigma}$ for $\sigma \ne \tau$ can only be nonzero if $h_{0}(\sigma)-h_{0}(\tau)$ is a positive integer, we conclude that $n_{\tau,\tau}=1$ are the only nonzero values and that $L_{0}(\tau)=M_{0}(\tau)$ for all $\tau$.

\subsection{WLOG $c>0$}  \label{c>0} There exists an isomorphism $$H_c(H_3,\mathfrak{h}) \to H_{-c}(H_3,\mathfrak{h}),$$ defined on the generators to be the identity on $\mathfrak{h}^*$ and $ \mathfrak{h}$ and to send $w \in H_3$ to ${\rm sign}_{H_3}(w)w$, where $ {\rm sign_{H_3}}:H_3 \to \left\{ -1,1\right\} $ is a representation of $H_3$, defined on any Coxeter group by sending the Coxeter generators to $-1$ (in our notation this is $\mathbf{1}_{-}$). Twisting by this isomorphism is an equivalence of categories between $\mathcal{O}_{c}$ and $\mathcal{O}_{-c}$. This equivalence exchanges  $M_{-c}(\mathbf{1}_{-} \otimes \tau)$ and $M_{c}(\tau)$, and consequently the same for their irreducible quotients $L_{c}(\tau)$. The map $\tau \to \mathbf{1}_{-}\times \tau$ is a permutation of $\hat{W}$. So, if we know the character formulas for all $L_{c}(\tau)$ in terms of $M_{c}(\sigma)$, just changing all the subscripts from $-$ to $+$ and from $+$ to $-$ will give us formulas for all $L_{-c}(\tau)$ in terms of $M_{-c}(\sigma)$.

\subsection{WLOG $c=r/d$, with $d=2,3,5,6,10$} \label{rationalc}
 The paper \cite{BE}, section 3.9, gathers results from \cite{DJO} and \cite{GGOR} to give exact conditions on $c\in \mathbb{C}$ such that the category $\mathcal{O}_c$ is not semisimple. Namely, it is shown that for $c>0$, $\mathcal{O}_c$ is not semisimple if and only if $c\in \mathbb{Q}$, and when written in a reduced form, its denominator is greater then $1$ and divides one of the degrees of basic invariants of the group $W$.

A module $M_{c}(\tau)$ is never a direct sum of two submodules, so if the category is semisimple, all $L_{c}(\tau)=M_{c}(\tau)$.

Degrees of basic invariants of $H_3$ are $2,6,10$; so all $c$ for which the character formulas are not trivial are reduced fractions of the form $r/d$, for $d=2,3,5,6,10$.

\subsection{WLOG $c=1/d$ with $d=3,5,6,10$, or $c=r/2$}\label{num=1,2}

\cite{GGOR} describes certain functors called $\mathrm{KZ}_{c}$ functors that associate to a representation from category $\mathcal{O}_{c}(W,\mathfrak{h})$ a representation of a certain Hecke algebra with parameter $e^{2\pi i c}$. For pairs of parameters $(c,c')$, these functors give rise to permutations $\varphi_{c,c'}$ of the set $\mathrm{Irrep} W$ of irreducible representations of $W$. We are interested in cases when $(c,c')=(1/d,r/d)$ with $d=3,5,6,10$ and $r$ relatively prime to $d$. The paper \cite{O} calculates these permutations explicitly for finite real reflection groups, and we can conclude from it that for $W=H_{3}$ and $r$ even, the permutation $\varphi=\varphi_{1/d,r/d}$ transposes the two 4-dimensional representations, i.e.  $\varphi(\mathbf{4}_{+})=\mathbf{4}_{-}$ and $\varphi(\mathbf{4}_{-})=\mathbf{4}_{+}$. The paper \cite{GG}, section 2.16, gives a formula for calculating $\varphi_{c,c'}$. It is explained there that one needs to consider the field extension of $\mathbb{Q}$ over which $\mathfrak{h}$ is defined (in our case, as $\mathfrak{h}\cong\mathbf{3}_{-}$, this is $\mathbb{Q}(\sqrt{5})$) and its Galois group (in our case this is $\mathbb{Z}_{2}$, as the only nontrivial field automorphism of $\mathbb{Q}(\sqrt{5})$ over $\mathbb{Q}$ is the one sending $\sqrt{5}$ to $-\sqrt{5}$). For a pair of parameters $(c,c')$ the map $e^{2\pi i c}\mapsto e^{2\pi i c'}$ induces an authomorphism of this field over $\mathbb{Q}$ and thus defines an element $g$ in the Galois group. The map $\varphi_{c,c'}$ then also permutes the irreducible representations of $W$ in a way that $g$ permutes their characters (this in our case corresponds to permuting the 3-dimensional representations, by transposing $\mathbf{3}_{-}$ and $\tilde{\mathbf{3}}_{-}$ and also transposing $\mathbf{3}_{+}$ and $\tilde{\mathbf{3}}_{+} $). In our case of $W=H_3$ and for pairs $(c,c')=(1/d,r/d)$ with $d=3,5,6,10$ and $r$ relatively prime to $d$, the only field extensions of $\mathbb{Q}$ by $e^{2\pi i r/d}$ that contain $\sqrt{5}$ are the in cases when $d=5$ or $d=10$. In those two cases, we calculate the effect of the maps $e^{2\pi i /d}\mapsto e^{2\pi i r/d}$ on $\sqrt{5}$ and get that the cases for which this map is a nontrivial element of the Galois group are $d=5$, and $r=2,3 \, (\mathrm{mod}\, 5)$, or when $d=10$, and $r=3,7 \, (\mathrm{mod}\, 10)$. These cases describe the permutation $\varphi_{1/d,r/d}$ completely.

The paper \cite{R1} then uses these KZ functors to give an equivalence of categories between $\mathcal{O}_{1/d}(W,\mathfrak{h})$ and $\mathcal{O}_{r/d}(W,\mathfrak{h})$ for above values of $r$ and $d$. As a result, describing the category $\mathcal{O}_{1/d}(W,\mathfrak{h})$ in those cases gives us results about $\mathcal{O}_{r/d}(W,\mathfrak{h})$.

The statement of Theorem 5.12 and Corollary 5.14 from \cite{R1}, after some necessary minor corrections, in our case of $W=H_{3}$, and using \cite{GG}, are as follows:

\begin{theorem}\label{equi1}
For $d=3,5,6,10$ and $r>0$ relatively prime with $d$, there exist equivalences of categories $$R_{r/d}:O_{1/d}\to O_{r/d}.$$ For each pair $(d,r)$ there exists a permutation $\varphi=\varphi_{1/d,r/d}$ such that  $$R_{r/d}(M_{1/d}(\tau))=M_{r/d}(\varphi(\tau)),$$ and$$R_{r/d}(L_{1/d}(\tau))=L_{r/d}(\varphi(\tau)).$$ 
Consequently, 
$$L_{1/d}(\tau)=\sum_{\sigma}n_{\tau,\sigma}M_{1/d}(\sigma) \Rightarrow L_{r/d}(\varphi(\tau))=\sum_{\sigma}n_{\tau,\sigma}M_{r/d}(\varphi(\sigma)).$$  $L_{r/d}(\varphi(\tau))$ is finite dimensional if and only if $L_{1/d}(\tau)$ is, with $\dim L_{r/d}(\varphi(\tau))=r^3\cdot \dim L_{1/d}(\tau)$. 

The permutation $\varphi=\varphi_{1/d,r/d}$ is given by:
\begin{itemize}
\item if $r=0 \, (\mathrm{mod}\, 2)$, $\varphi(\mathbf{4}_{-})=\mathbf{4}_{+}$, $\varphi(\mathbf{4}_{+})=\mathbf{4}_{-}$
\item if $d=5$ and $r=2,3 \, (\mathrm{mod}\, 5)$, $\varphi(\mathbf{3}_{-})=\tilde{\mathbf{3}}_{-}$, $\varphi(\tilde{\mathbf{3}}_{-})=\mathbf{3}_{-}$, $\varphi(\mathbf{3}_{+})=\tilde{\mathbf{3}}_{+}$, $\varphi(\tilde{\mathbf{3}}_{+})=\mathbf{3}_{+}$
\item if $d=10$ and $r=3,7 \, (\mathrm{mod}\, 10)$, $\varphi(\mathbf{3}_{-})=\tilde{\mathbf{3}}_{-}$, $\varphi(\tilde{\mathbf{3}}_{-})=\mathbf{3}_{-}$, $\varphi(\mathbf{3}_{+})=\tilde{\mathbf{3}}_{+}$, $\varphi(\tilde{\mathbf{3}}_{+})=\mathbf{3}_{+}$
\item $\varphi(\tau)=\tau$ for modules and indices not listed above.
\end{itemize}

\end{theorem}

\subsection{WLOG $c=1/d$ with $d=2,3,5,6,10$ or $c=r/2$, for finitely many $r$}\label{num=1}
Even though the functors from \cite{R1} are proved to be equivalences of categories only for $c=r/d$, $d\ne 2$, they are conjectured to be equivalences for $d=2$ as well. In absence of the proof of that, we can use other functors for reducing the number of $c=r/2$ parameters we need to check. For this, we will recall the results of \cite{BEG1}, and the notion of translation functors. 

Denote $H_{c}=H_{c}(W,\mathfrak{h})$, for $c$ a constant.

Let $e_{+}=\frac{1}{|W|}\sum_{w\in W} w\in \mathbb{C}W$ and $e_{-}=\frac{1}{|W|}\sum_{w\in W} \mathrm{sign}_{W}(w)w\in \mathbb{C}W$ be projections to $W$ invariants and anti-invariants, respectively. As $\mathbb{C}W$ is a subalgebra of $H_{c}$ for every $c$, we can consider $e_{+}$ and $e_{-}$ as elements of any $H_{c}$. The paper \cite{BEG1} then shows there is an isomorphism of filtered algebras $$\phi_{c}:e_{+}H_{c}e_{+}\to e_{-}H_{c+1}e_{-}.$$ It induces the natural equivalence between categories of their representations $$\Phi_{c}:\mathcal{O}(e_{+}H_{c}e_{+})\to \mathcal{O}(e_{-}H_{c+1}e_{-}).$$

Next, for $\varepsilon\in \{ +,-\}$, one defines functors $$F_{c}^{\varepsilon}:\mathcal{O}_{c}\to \mathcal{O}_{c}(e_{\varepsilon}H_{c}e_{\varepsilon}) \qquad G_{c}^{\varepsilon}:\mathcal{O}_{c}\to \mathcal{O}_{c}(e_{\varepsilon}H_{c}e_{\varepsilon})$$ by $$F_{c}^{\varepsilon}(V)=e_{\varepsilon}V \qquad G_{c}^{\varepsilon}(Y)=H_{c}e_{\varepsilon}\otimes_{e_{\varepsilon}H_{c}e_{\varepsilon}}Y.$$

Denote by $\mathcal{O}_{c}^{\varepsilon}$ the full subcategory of $\mathcal{O}_{c}$ consisting of modules such that $e_{\varepsilon}V=0$. This is a Serre subcategory so the quotients  $\mathcal{O}_{c}/\mathcal{O}_{c}^{\varepsilon}$ make sense. Using the fact that for $V$ a simple module in $\mathcal{O}_{c}$, either $e_{\varepsilon}V=$ or $H_{c}e_{\varepsilon}V=V$, it is easy to see that $F_{c}^{\varepsilon}$ and $G_{c}^{\varepsilon}$, now understood as functors to and from quotient categories, are mutually inverse equivalences of categories between $\mathcal{O}_{c}/\mathcal{O}_{c}^{\varepsilon}$ and $\mathcal{O}(e_{\varepsilon}H_{c}e_{\varepsilon})$. So, we have the following diagram of equivalences of categories:
$$\begin{CD} 
\mathcal{O}_{c}/\mathcal{O}_{c}^+ @>F_{c}^{+}>> \mathcal{O}(e_{+}H_{c}e_{+}) @>\Phi_{c}>> \mathcal{O}(e_{-}H_{c+1}e_{-})  @>G_{c+1}^{-}>> \mathcal{O}_{c+1}/\mathcal{O}_{c+1}^-
\end{CD} $$

Consider the compositions of these functors $S_{c}^+=G_{c+1}^{-}\circ \Phi_{c}\circ F_{c}^+$, and $S_{c}^-=G_{c}^{+}\circ \Phi_{c}^{-1} \circ F_{c+1}^-$. These are mutually inverse equivalences of categories between $\mathcal{O}_{c}/\mathcal{O}_{c}^+$ and $\mathcal{O}_{c+1}/\mathcal{O}_{c+1}^-$. Moreover, if $\mathcal{O}_{c}^+=0$, then $S_{c}^+$ is an equivalence of categories between $\mathcal{O}_{c}$ and $\mathcal{O}_{c+1}$; because $\mathcal{O}_{c}$ and $\mathcal{O}_{c+1}$ have the same number of simple objects (equal to the number of irreducible representations of $W$), and if $\mathcal{O}_{c+1}^-\ne 0$, then $\mathcal{O}_{c+1}/\mathcal{O}_{c+1}^-$ has strictly less. Due to \cite{BE} and \cite{Chm}, all the infinite dimensional $L_{c}(\tau)$ contain a copy of a trivial representation, i.e. $e_{+}L_{c}(\tau)\ne 0$. 

From \cite{BEG1}, \cite{GG} and \cite{O} we can also conclude what $S_{c}^+$ does to standard and irreducible modules:
\begin{theorem} \label{equi2}
For $c=r/2,$ $r$ odd, $r>0$, the functor $$S_{c}^{+}:\mathcal{O}_{c}/\mathcal{O}_{c}^+\to \mathcal{O}_{c+1}/\mathcal{O}_{c+1}^-$$ is an equivalence of categories, with 
$$S_{c}^{+}\overline{M_{c}(\tau)}= \overline {M_{c+1}(\varphi(\tau))} \qquad S_{c}^{+} \overline {L_{c}(\tau)}= \overline {L_{c+1}(\varphi(\tau))}, $$ where $\varphi$ is a permutation of irreducible representations of $W$ from the previous subsection (in case $W=H_{3}$, $\varphi (\mathbf{4}_{+})=\mathbf{4}_{-},\varphi (\mathbf{4}_{-})=\mathbf{4}_{+}$, and $\varphi (\tau)=\tau$ for all other $\tau$), and $\overline {M_{c}(\tau)}, \overline {L_{c}(\tau)}$ denote the images of $M_{c}(\tau)$ and $L_{c}(\tau)$ in the quotient of categories $\mathcal{O}_{c}/\mathcal{O}_{c}^+$ and $\mathcal{O}_{c+1}/\mathcal{O}_{c+1}^-$.

If $\mathcal{O}_{c}^+=0$, then $\mathcal{O}_{c+1}^-=0$ and $S_{c}^{+}$ is an equivalence $\mathcal{O}_{c}\to \mathcal{O}_{c+1}$. This happens for $c$ large enough.
\end{theorem}

\section{Calculations for $c=1/10$}\label{ch1/10}

\begin{theorem}
 Irreducible representations in category $\mathcal{O}_{1/10}(H_3,\mathfrak{h})$ have the following descriptions in the Grothendieck group:
\begin{eqnarray*}
L(\mathbf{1}_{+}) & = & M(\mathbf{1}_{+})- M(\mathbf{3}_{-})+M(\mathbf{3}_{+})-M(\mathbf{1}_{-}) \\
L(\mathbf{1}_{-}) & = & M(\mathbf{1}_{-}) \\
L(\mathbf{3}_{+})& = &M(\mathbf{3}_{+})-M(\mathbf{1}_{-})\\
L(\mathbf{3}_{-})& = &M(\mathbf{3}_{-})-M(\mathbf{3}_{+})+M(\mathbf{1}_{-}) \\
L(\widetilde{\mathbf{3}}_{+})& = & M(\widetilde{\mathbf{3}}_{+})\\
L(\widetilde{\mathbf{3}}_{-})& = & M( \widetilde{\mathbf{3}}_{-})\\
L(\mathbf{4}_{+})& = &M(\mathbf{4}_{+})\\
L(\mathbf{4}_{-})& = &M(\mathbf{4}_{-})\\
L(\mathbf{5}_{+})& = &M(\mathbf{5}_{+})\\
L(\mathbf{5}_{-})& = &M(\mathbf{5}_{-})
\end{eqnarray*}
Among these representations only $L(\mathbf{1}_{+})$ is finite dimensional, with $\mathrm{ch}_{L(\mathbf{1}_{+})}(w,t)=1.$
\label{1/10}
\end{theorem}

The rest of this chapter is the proof of this theorem. Let us first calculate the constants $h_{1/10}(\tau)=\frac{3}{2}-\frac{1}{10}\sum_{s\in S} s|_{\tau}$ (see Table \ref{h1/10}).
\begin{table}[h]
\begin{center}
\begin{tabular}{|c|c|c|c|c|c|c|c|c|c|} \hline
$\mathbf{1}_{+}$ & $\mathbf{1}_{-}$ & $\mathbf{3}_{+}$ & $\mathbf{3}_{-}$ & $\widetilde{\mathbf{3}}_+$ & $\widetilde{\mathbf{3}}_-$ & $\mathbf{4}_{+}$ & $\mathbf{4}_{-}$ & $\mathbf{5}_{+}$ & $\mathbf{5}_{-}$ \\ \hline
0 & 3 & 2 & 1 & 2 & 1 & 3/2 & 3/2 & 6/5 & 9/5 \\ \hline
\end{tabular}
\end{center}
\caption{$h_{1/10}(\tau)$}
\label{h1/10}
\end{table}

Using lemma \ref{even} we immediately conclude: 
$$L(\mathbf{4}_{+})=M(\mathbf{4}_{+})$$
$$L(\mathbf{4}_{-})=M(\mathbf{4}_{-})$$
$$L(\mathbf{5}_{+})=M(\mathbf{5}_{+})$$
$$L(\mathbf{5}_{-})=M(\mathbf{5}_{-})$$

Mark the lowest weights of other modules on the real line as 

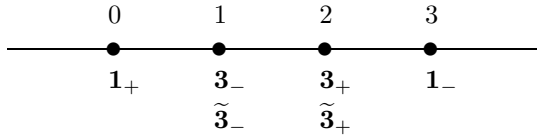
\begin{center}
\begin{picture}(200,50)
\put(0,30){\line(1,0){200}}
\put(40,30){\circle*{5}}
\put(80,30){\circle*{5}}
\put(120,30){\circle*{5}}
\put(160,30){\circle*{5}}
\put(38,40){$0$}
\put(78,40){$1$}
\put(118,40){$2$}
\put(158,40){$3$}
\put(38,15){$\mathbf{1}_{+}$}
\put(78,15){$\mathbf{3}_{-}$}
\put(78,0){$\widetilde{\mathbf{3}}_-$}
\put(118,15){$\mathbf{3}_{+}$}
\put(118,0){$\widetilde{\mathbf{3}}_+$}
\put(158,15){$\mathbf{1}_{-}$}
\end{picture}
\end{center}

This picture represents Lemma \ref{even} graphically, meaning that $n_{\tau,\sigma}$ can be nonzero only if both $\tau$ and $\sigma$ are represented on the line, with $\sigma$ to the right of $\tau$.  From this we can also immediately conclude that $$L(\mathbf{1}_{-})=M(\mathbf{1}_{-}).$$

To calculate character formulas for $L(\mathbf{3}_{+})$ and $L(\widetilde{\mathbf{3}}_{+})$, we will use Lemma \ref{Lemma35}. First calculate the decomposition into irreducible $H_3$ representations of $\mathfrak{h}^* \otimes \mathbf{3}_{+}$ and $\mathfrak{h}^* \otimes\widetilde{\mathbf{3}}_{+}$. A computation with characters of finite group $H_3$ (see \cite{FH} and Table \ref{H3CharTable}) gives 
$$\mathfrak{h}^* \otimes \mathbf{3}_{+}= \mathbf{3}_{-}\otimes \mathbf{3}_{+}\cong 
\mathbf{1}_{-}\oplus \mathbf{5}_{-}\oplus \mathbf{3}_{-}
$$
and 
$$\mathfrak{h}^* \otimes\widetilde{\mathbf{3}}_{+}= \mathbf{3}_{-}\otimes \widetilde{\mathbf{3}}_{+}\cong  \mathbf{4}_{-}\oplus  \mathbf{5}_{-} .$$

Lemma \ref{Lemma35} now implies that the subrepresentation $\sigma=\mathbf{1}_{-}\subseteq \mathfrak{h}^* \otimes \mathbf{3}_{+}$ consists of singular vectors, and hence that $M(\mathbf{1}_{-})$ is a subrepresentation of $M(\mathbf{3}_{+})$. It is the maximal proper subrepresentation, and it follows that $L(\mathbf{3}_{+})\cong M(\mathbf{3}_{+}) / M(\mathbf{1}_{-})$, so in the Grothendieck group
 $$L(\mathbf{3}_{+})=M(\mathbf{3}_{+})-M(\mathbf{1}_{-}).$$

On the other hand, decomposition of $\mathfrak{h}^* \otimes\widetilde{\mathbf{3}}_{+}$ doesn't have $\mathbf{1}_{-}$ as a subrepresentation, so $$L(\widetilde{\mathbf{3}}_{+})=M(\widetilde{\mathbf{3}}_{+}).$$

Next, using the decomposition $\mathrm{S}^2\mathfrak{h}^*\cong \mathbf{1}_{+}\oplus \mathbf{5}_{+}$, let us decompose two more $H_3$ representations:
$$\mathfrak{h}^*\otimes \widetilde{\mathbf{3}}_{-}\cong \mathbf{4}_{+}\oplus \mathbf{5}_{+}$$
$$\mathrm{S}^2\mathfrak{h}^*\otimes  \widetilde{\mathbf{3}}_{-}\cong  (\mathbf{1}_{+}\oplus \mathbf{5}_{+})\otimes \widetilde{\mathbf{3}}_{-} \cong 2\widetilde{\mathbf{3}}_{-}  \oplus \mathbf{3}_{-}\oplus \mathbf{4}_{-}\oplus \mathbf{5}_{-}.$$

Because neither $\mathbf{3}_{+} $ nor $\widetilde{\mathbf{3}}_{+} $ appear as subrepresentations of $\mathfrak{h}^*\otimes \widetilde{\mathbf{3}}_{-}$, nor does $\mathbf{1}_{-}$ appear in the decomposition of $\mathrm{S}^2\mathfrak{h}^*\otimes  \widetilde{\mathbf{3}}_{-}$, the module $M( \widetilde{\mathbf{3}}_{-})$ must be simple: $$L( \widetilde{\mathbf{3}}_{-})=  M( \widetilde{\mathbf{3}}_{-}).$$

Corresponding decompositions for $\mathbf{3}_{-}$ are
$$\mathfrak{h}^*\otimes \mathbf{3}_{-}\cong \mathbf{1}_{+}\oplus \mathbf{3}_{+}\oplus \mathbf{5}_{+}$$
$$\mathrm{S}^2\mathfrak{h}^*\otimes \mathbf{3}_{-}\cong 2\cdot \mathbf{3}_{-}\oplus  \widetilde{\mathbf{3}}_{-} \oplus \mathbf{4}_{-}\oplus \mathbf{5}_{-}.$$

From the first of these formulas and using lemma \ref{Lemma35} we can now conclude that $\mathbf{3}_{+} \subseteq M(\mathbf{3}_{-})[2]$ consists of singular vectors, so it generates a $H(H_3, \mathfrak{h})$ subrepresentation. $\mathbf{1}_{-}$ doesn't appear in the decomposition of $\mathrm{S}^2\mathfrak{h}^*\otimes \mathbf{3}_{-}$, so the subrepresentation  generated by $\mathbf{3}_{+}$ is the whole $J(\mathbf{3}_{-})$. Looking at the computations for $L(\mathbf{3}_{+})$ we see that the only lowest weight representations with lowest weight $\mathbf{3}_{+}$ are $M(\mathbf{3}_{+})$ and $L(\mathbf{3}_{+})=M(\mathbf{3}_{+})-M(\mathbf{1}_{-})$. Thus in Grothendieck group, $L(\mathbf{3}_{-})=M(\mathbf{3}_{-})-M(\mathbf{3}_{+})+n_{\mathbf{3}_{-},\mathbf{1}_{-}}M(\mathbf{1}_{-})$, for $n_{\mathbf{3}_{-},\mathbf{1}_{-}}=0$ or $n_{\mathbf{3}_{-},\mathbf{1}_{-}}=1$. To see which one of these it is, notice that $\mathbf{1}_{-}$ doesn't appear as an $H_3$ subrepresentation in $M(\mathbf{3}_{-})[3]\cong \mathrm{S}^2\mathfrak{h}^*\otimes \mathbf{3}_{-}$, but it does in $M(\mathbf{3}_{+})[3]\cong \mathrm{S}^1\mathfrak{h}^*\otimes \mathbf{3}_{+}$. That means that $M(\mathbf{3}_{+})$ cannot be a submodule of $M(\mathbf{3}_{-})$; so $J(\mathbf{3}_{-})=L(\mathbf{3}_{+})$, and the expression in Grothendieck group we wanted is 
$$L(\mathbf{3}_{-})=M(\mathbf{3}_{-})-M(\mathbf{3}_{+})+M(\mathbf{1}_{-}).$$

None of the modules considered so far in this chapter is finite dimensional. An easy way to see that is to  consider them as $\mathfrak{sl}_2$ representations. The lowest occuring weights are then given by Table \ref{h1/10}. Every finite dimensional $\mathfrak{sl}_2$ representation will have integral weights, with the lowest one being less or equal then $0$. As none of the lowest weights of these modules is a nonpositive integer, they are not finite dimensional. 

This doesn't apply to the one module we still didn't describe, $L(\mathbf{1}_{+})$. Its lowest $\mathbf{h}$ weight is $0$, so it could be finite dimensional in case it was a trivial one dimensional module. That is exactly what happens: it is easy to see by direct calculation that setting $x=0, y=0, w=1, x\in \mathfrak{h}^*, y\in \mathfrak{h}, w\in H_3$ defines an action of $H(H_3,\mathfrak{h})$ on $\mathbb{C}$. So, there is a trivial module at $c=1/10$, whose lowest weight is $\mathbf{1}_{+}$, and it has to be $L(\mathbf{1}_{+})$. This computation appears in \cite{BEG2}, Prop 2.1. 

The character of $L(\mathbf{1}_{+})$ is naturally $1$; to express it in terms of characters of $M(\sigma)$, we count the dimensions of $\mathbf{h}$ weight spaces. Clearly the copy of $\mathbf{3}_{-}\subseteq M(\mathbf{1}_{+})[1]$ consists of singular vectors, spanning either $M(\mathbf{3}_{-})$ or  $L(\mathbf{3}_{-})$. To see which one it is, look at the next $\mathbf{h}$ weight space, where $\dim M(\mathbf{3}_{-})[2] =9>6= \dim M(\mathbf{1}_{+})[2]$. So, the submodule with the lowest weight in $\mathbf{h}$ weight space $1$ is $L(\mathbf{3}_{-})$. The dimensions of all higher $\mathbf{h}$ weight spaces of $M(\mathbf{1}_{+})/L(\mathbf{3}_{-})$ are $0$, so $J(\mathbf{1}_{+})=L(\mathbf{3}_{-})$ and 
 $$L(\mathbf{1}_{+})= M(\mathbf{1}_{+})- L(\mathbf{3}_{-})= M(\mathbf{1}_{+})- M(\mathbf{3}_{-})+M(\mathbf{3}_{+})-M(\mathbf{1}_{-}).$$

\section{Calculations for $c=1/6$}\label{ch1/6}

\begin{theorem}
 Irreducible representations in category $\mathcal{O}_{1/6}(H_3,\mathfrak{h})$ have the following descriptions in the Grothendieck group:
\begin{eqnarray*}
L(\mathbf{1}_{+})&=&M(\mathbf{1}_{+})-M(\mathbf{5}_{+})+M(\mathbf{5}_{-})-M(\mathbf{1}_{-})\\
L(\mathbf{1}_{-})&=&M(\mathbf{1}_{-})\\
L(\mathbf{3}_{+})&=&M(\mathbf{3}_{+})\\
L(\mathbf{3}_{-})&=&M(\mathbf{3}_{-})\\
L(\widetilde{\mathbf{3}}_+)&=&M(\widetilde{\mathbf{3}}_+)\\
L(\widetilde{\mathbf{3}}_-)&=&M(\widetilde{\mathbf{3}}_-)\\
L(\mathbf{4}_{+})&=&M(\mathbf{4}_{+})\\
L(\mathbf{4}_{-})&=&M(\mathbf{4}_{-})\\
L(\mathbf{5}_{+})&=&M(\mathbf{5}_{+})-M(\mathbf{5}_{-})+M(\mathbf{1}_{-})\\
L(\mathbf{5}_{-})&=&M(\mathbf{5}_{-})-M(\mathbf{1}_{-})
\end{eqnarray*}
Only $L(\mathbf{1}_{+})$ among these representations is finite dimensional, with character $\mathrm{ch}_{L(\mathbf{1}_{+})}=\chi_{\mathbf{1}_{+}}t^{-1} + \chi_{\mathbf{3}_{-}}+\chi_{\mathbf{1}_{+}}t$.
\label{1/6}
\end{theorem}

Let us again first calculate the constants $h_{1/6}(\tau)=\frac{3}{2}-\frac{1}{6}\sum_{s\in S} s|_{\tau}$ (see Table \ref{h1/6}).
\begin{table}[h]
\begin{center}
\begin{tabular}{|c|c|c|c|c|c|c|c|c|c|} \hline
$\mathbf{1}_{+}$ & $\mathbf{1}_{-}$ & $\mathbf{3}_{+}$ & $\mathbf{3}_{-}$ & $\widetilde{\mathbf{3}}_+$ & $\widetilde{\mathbf{3}}_-$ & $\mathbf{4}_{+}$ & $\mathbf{4}_{-}$ & $\mathbf{5}_{+}$ & $\mathbf{5}_{-}$ \\ \hline
-1 & 4 & 7/3 & 2/3 & 7/3 & 2/3 & 3/2 & 3/2 & 1 & 2 \\ \hline
\end{tabular}
\end{center}
\caption{$h_{1/6}(\tau)$}
\label{h1/6}
\end{table}

We immediately conclude that $M(\mathbf{3}_{+})$, $M(\mathbf{3}_{-})$, $M(\widetilde{\mathbf{3}}_+)$, $M(\widetilde{\mathbf{3}}_-)$, $M(\mathbf{4}_{+})$ and $M(\mathbf{4}_{-})$ are simple.

The remaining modules have lowest weights represented in the following picture:
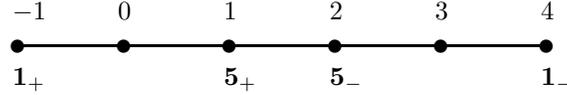
\begin{center}
\begin{picture}(200,50)
\put(0,30){\line(1,0){200}}
\put(0,30){\circle*{5}}
\put(40,30){\circle*{5}}
\put(80,30){\circle*{5}}
\put(120,30){\circle*{5}}
\put(160,30){\circle*{5}}
\put(200,30){\circle*{5}}
\put(-2,40){$-1$}
\put(38,40){$0$}
\put(78,40){$1$}
\put(118,40){$2$}
\put(158,40){$3$}
\put(198,40){$4$}
\put(-2,15){$\mathbf{1}_{+}$}
\put(78,15){$\mathbf{5}_{+}$}
\put(118,15){$\mathbf{5}_{-}$}
\put(198,15){$\mathbf{1}_{-}$}
\end{picture}
\end{center}
 
 So, $M(\mathbf{1}_{-})$ is also simple. 

Calculate $$\mathrm{S}^2\mathfrak{h}^*\otimes \mathbf{5}_{-}\cong \mathbf{1}_{-} \oplus \mathbf{3}_{-}\oplus  \widetilde{\mathbf{3}}_{-}  \oplus  2  \mathbf{4}_{-} \oplus 3 \mathbf{5}_{-}$$
From this we can conclude that $L(\mathbf{5}_{-})=M(\mathbf{5}_{-})-n_{\mathbf{5}_{-},\mathbf{1}_{-}}\cdot M(\mathbf{1}_{-})$, with $n_{\mathbf{5}_{-},\mathbf{1}_{-}}\in \lbrace 0,1 \rbrace$. It is possible to deduce $n_{\mathbf{5}_{-},\mathbf{1}_{-}}$ from the rank of the contravariant form $B$ restricted to $M(\mathbf{5}_{-})[4]$, but we will use a less direct argument here.

Let us focus on $\mathbf{5}_{+}$ for a while. We notice that $$\mathfrak{h}^*\otimes \mathbf{5}_{+}\cong \mathbf{3}_{-} \oplus \widetilde{\mathbf{3}}_{-} \oplus \mathbf{4}_{-}\oplus \mathbf{5}_{-},$$ so by Lemma \ref{Lemma35}, the $H_3$ subrepresentation $\mathbf{5}_{-}$ consists of lowest weight vectors. We know from the previous paragraph that if $n_{\mathbf{5}_{-},\mathbf{1}_{-}}=0$, then there is just one representation with lowest weight $\mathbf{5}_{-}$, that is $M(\mathbf{5}_{-})$, and if $n_{\mathbf{5}_{-},\mathbf{1}_{-}}=1$  there are two, namely the standard one $M(\mathbf{5}_{-})$ and the irreducible one $M(\mathbf{5}_{-})-M(\mathbf{1}_{-})$. The module $M(\mathbf{5}_{+})$ can also have a $b$-dimensional space of singular vectors in $M(\mathbf{5}_{+})[4]$, $b\in \mathbb{N}_0$. So, the expression for $L(\mathbf{5}_{+})$ is either $$L(\mathbf{5}_{+})=M(\mathbf{5}_{+})-M(\mathbf{5}_{-})-bM(\mathbf{1}_{-}),$$ or $$L(\mathbf{5}_{+})=M(\mathbf{5}_{+})-M(\mathbf{5}_{-})+(n_{\mathbf{5}_{-},\mathbf{1}_{-}}-b)M(\mathbf{1}_{-}).$$ Now use the decompositions 
$$\mathrm{S}^3\mathfrak{h}^*\otimes \mathbf{5}_{+}\cong 3\cdot \mathbf{3}_{-}\oplus 3 \cdot \widetilde {\mathbf{3}}_{-}\oplus 3 \cdot \mathbf{4}_{-}\oplus 4\cdot \mathbf{5}_{-}$$
$$\mathrm{S}^2\mathfrak{h}^*\otimes \mathbf{5}_{-}\cong \mathbf{1}_{-}\oplus \mathbf{3}_{-}\oplus \widetilde {\mathbf{3}}_{-}\oplus 2\cdot \mathbf{4}_{-}\oplus 3\cdot \mathbf{5}_{-}$$ to deduce that the graded piece of $L(\mathbf{5}_{+})$ with $\mathbf{h}$ weight 4 has, in the Grothendieck group of $H_3$ representations, one of the following two decompositions:
\begin{eqnarray*} L(\mathbf{5}_{+})[4]& = &M(\mathbf{5}_{+})[4]-M(\mathbf{5}_{-})[4]-b\cdot M(\mathbf{1}_{-})[4] \\
& = & \mathrm{S}^3\mathfrak{h}^*\otimes \mathbf{5}_{+} -  \mathrm{S}^2\mathfrak{h}^*\otimes \mathbf{5}_{-} -b\cdot  \mathbf{1}_{-} \\
& = & 3\cdot \mathbf{3}_{-} + 3\cdot \widetilde {\mathbf{3}}_{-} + 3\cdot \mathbf{4}_{-}+ 4\cdot \mathbf{5}_{-} - 
 \mathbf{1}_{-}- \mathbf{3}_{-} - \widetilde {\mathbf{3}}_{-}- 2\cdot \mathbf{4}_{-}- 3\cdot \mathbf{5}_{-} -b\cdot  \mathbf{1}_{-} \\
& = & (-1-b)\cdot \mathbf{1}_{-} + 2\cdot \mathbf{3}_{-} + 2\cdot \widetilde {\mathbf{3}}_{-} + \mathbf{4}_{-}+ \mathbf{5}_{-},\\
\end{eqnarray*}
or
\begin{eqnarray*} L(\mathbf{5}_{+})[4]& = &M(\mathbf{5}_{+})[4]-M(\mathbf{5}_{-})[4]+(a-b)\cdot M(\mathbf{1}_{-})[4] \\
& = & (n_{\mathbf{5}_{-},\mathbf{1}_{-}}-1-b)\cdot \mathbf{1}_{-} + 2\cdot \mathbf{3}_{-} + 2\cdot \widetilde {\mathbf{3}}_{-} + \mathbf{4}_{-}+ \mathbf{5}_{-}.\\
\end{eqnarray*}
These are decompositions of an actual representation of $H_3$, so all the coefficients need to be nonnegative integers. $-1-b$ can never be more than $-1$, so the correct decomposition is the second one, $n_{\mathbf{5}_{-},\mathbf{1}_{-}}=1,b=0$, and the correct formulas for both irreducible modules are
$$L(\mathbf{5}_{-})=M(\mathbf{5}_{-})-M(\mathbf{1}_{-})$$
$$L(\mathbf{5}_{+})=M(\mathbf{5}_{+})-M(\mathbf{5}_{-})+M(\mathbf{1}_{-}).$$

To describe the module $L(\mathbf{1}_{+})$ we use Theorem \ref{support}. It says that its support, when viewed as a $\mathbb{C}[\mathfrak{h}]$ module, is the set of $a\in\mathfrak{h}$ such that  $\# \lbrace i | 6 \text{ divides } d_i(W)  \rbrace= \# \lbrace i | 6 \text{ divides } d_i(W_a)  \rbrace$. Degrees of $H_3$ are $2,6,10$, so the size of that set is $1$. Maximal parabolic subgroups $W_a$ of $H_3$ are Coxeter groups obtained by deleting a node from the Coxeter graph of $H_3$, so they are $A_1\times A_1$, $I_{2}(5)$ and $A_2$. The degrees of their basic invariants are: $d_1(A_1)=2, d_1(I_2(5))=2, d_2(I_2(5))=5, d_1(A_2)=2, d_2(A_2)=3$. Since $6$ doesn't divide any of them, theorem implies that support of $L(\mathbf{1}_{+})$ is just $0\in\mathfrak{h}$. Thus, its Hilbert-Poincar\'e series doesn't have a pole at t=1, so it is a polynomial, and $L(\mathbf{1}_{+})$ is finite dimensional.

We know that the expression $L(\mathbf{1}_{+})$ in the Grothendieck group is of the form $$L(\mathbf{1}_{+})=M(\mathbf{1}_{+})+n_{\mathbf{1}_{+},\mathbf{5}_{+}} M(\mathbf{5}_{+})+n_{\mathbf{1}_{+},\mathbf{5}_{-}} M(\mathbf{5}_{-})+n_{\mathbf{1}_{+},\mathbf{1}_{-}} M(\mathbf{1}_{-}).$$ The characters of these representations relate in the same way. Substituting the known expression for character of $M_{c}(\tau)$ and evaluating at $w=1$, we get that
 $$\mathrm{ch}_{L(\mathbf{1}_{+})}(\mathrm{Id},t)=\frac{t^{-1}}{(1-t)^3}+n_{\mathbf{1}_{+},\mathbf{5}_{+}}\cdot \frac{5t}{(1-t)^3}+n_{\mathbf{1}_{+},\mathbf{5}_{-}}\cdot \frac{t^{2}}{(1-t)^3}+n_{\mathbf{1}_{+},\mathbf{1}_{-}}\cdot \frac{t^{4}}{(1-t)^3}$$
must be regular at $t=1$, i.e. that 
 $$t^{-1}+n_{\mathbf{1}_{+},\mathbf{5}_{+}}\cdot 5t+n_{\mathbf{1}_{+},\mathbf{5}_{-}}\cdot t^{2}+n_{\mathbf{1}_{+},\mathbf{1}_{-}}\cdot t^{4}$$
must vanish to order $3$ at $t=1$. Solving this system we get that the only case when this happens is $n_{\mathbf{1}_{+},\mathbf{5}_{+}}=n_{\mathbf{1}_{+},\mathbf{1}_{-}}=-1, n_{\mathbf{1}_{+},\mathbf{5}_{-}}=1$, so the Grothendieck group expression is, as claimed,
$$L(\mathbf{1}_{+})=M(\mathbf{1}_{+})-M(\mathbf{5}_{+})+M(\mathbf{5}_{-})-M(\mathbf{1}_{-}).$$
It is now an easy computation of $H_3$ characters in $\mathrm{S}^2\mathfrak{h}^*\otimes \mathbf{1}_{+}$ to see that this is also equal to $$\chi_{\mathbf{1}_{+}}t^{-1} + \chi_{\mathbf{3}_{-}}+\chi_{\mathbf{1}_{+}}t.$$

Looking at lowest $\mathbf{h}$ weights again, we see that no module other then $L(\mathbf{1}_{+})$ can be finite dimensional, which completes the proof.

\section{Calculations for $c=1/5$}\label{ch1/5}

\begin{theorem}
 Irreducible representations in category $\mathcal{O}_{1/5}(H_3,\mathfrak{h})$ have the following descriptions in the Grothendieck group:
\begin{eqnarray*}
L( \mathbf{1}_{+})&=&M( \mathbf{1}_{+})-M( \mathbf{4}_{-})+M( \widetilde{\mathbf{3}}_{+})\\
L(\mathbf{1}_{-})&=&M(\mathbf{1}_{-})\\
L(\mathbf{3}_{+})&=&M(\mathbf{3}_{+})\\
L(\mathbf{3}_{-})&=&M(\mathbf{3}_{-})\\
L(\widetilde{\mathbf{3}}_{+})&=&M(\widetilde{\mathbf{3}}_{+})\\
L(\widetilde{\mathbf{3}}_{-})&=&M(\widetilde{\mathbf{3}}_{-})-M(\mathbf{4}_{+})+M(\mathbf{1}_{-})\\
L(\mathbf{4}_{+})&=&M(\mathbf{4}_{+})-M(\mathbf{1}_{-})\\
L(\mathbf{4}_{-})&=&M(\mathbf{4}_{-})-M(\widetilde{\mathbf{3}}_{+})\\
L(\mathbf{5}_{+})&=&M(\mathbf{5}_{+})\\
L(\mathbf{5}_{-})&=&M(\mathbf{5}_{-})
\end{eqnarray*}
None of these representations is finite dimensional.
\label{1/5}
\end{theorem}

In this case, $h_{1/5}(\tau)=\frac{3}{2}-\frac{1}{5}\sum_{s\in S} s|_{\tau}$ are as follows (see Table \ref{h1/5}):
\begin{table}[h!]
\begin{center}
\begin{tabular}{|c|c|c|c|c|c|c|c|c|c|} \hline
$\mathbf{1}_{+}$ & $\mathbf{1}_{-}$ & $\mathbf{3}_{+}$ & $\mathbf{3}_{-}$ & $\widetilde{\mathbf{3}}_+$ & $\widetilde{\mathbf{3}}_-$ & $\mathbf{4}_{+}$ & $\mathbf{4}_{-}$ & $\mathbf{5}_{+}$ & $\mathbf{5}_{-}$ \\ \hline
-3/2 & 9/2 & 5/2 & 1/2 & 5/2 & 1/2 & 3/2 & 3/2 & 9/10 & 21/10 \\ \hline
\end{tabular}
\end{center}
\caption{$h_{1/5}(\tau)$}
\label{h1/5}
\end{table}

An observation we can immediately make by restricting the representations to the $\mathfrak{sl}_2$ subalgebra is that there are no finite dimensional modules at $c=1/5$, because those would have integral weights. We can also immediately say that $M(\mathbf{5}_{+})$ and $M(\mathbf{5}_{-})$ are simple.

Taking into consideration Lemma \ref{even}, draw the remaining $8$ representations schematically as
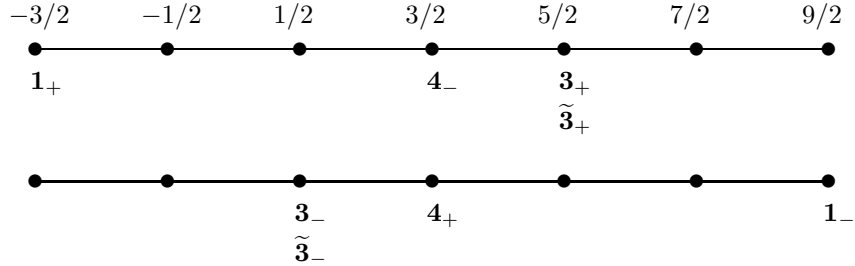
\begin{center}
\begin{picture}(300,105)
\put(0,80){\line(1,0){300}}
\put(0,30){\line(1,0){300}}
\put(0,30){\circle*{5}}
\put(50,30){\circle*{5}}
\put(100,30){\circle*{5}}
\put(150,30){\circle*{5}}
\put(200,30){\circle*{5}}
\put(250,30){\circle*{5}}
\put(300,30){\circle*{5}}
\put(0,80){\circle*{5}}
\put(50,80){\circle*{5}}
\put(100,80){\circle*{5}}
\put(150,80){\circle*{5}}
\put(200,80){\circle*{5}}
\put(250,80){\circle*{5}}
\put(300,80){\circle*{5}}
\put(-10,90){$-3/2$}
\put(40,90){$-1/2$}
\put(90,90){$1/2$}
\put(140,90){$3/2$}
\put(190,90){$5/2$}
\put(240,90){$7/2$}
\put(290,90){$9/2$}
\put(-2,65){$\mathbf{1}_{+}$}
\put(298,15){$\mathbf{1}_{-}$}
\put(198,65){$\mathbf{3}_{+}$}
\put(198,50){$\widetilde{\mathbf{3}}_{+}$}
\put(98,15){$\mathbf{3}_{-}$}
\put(98,0){$\widetilde{\mathbf{3}}_{-}$}
\put(148,65){$\mathbf{4}_{-}$}
\put(148,15){$\mathbf{4}_{+}$}
\end{picture}
\end{center}
 
This picture means that $n_{\tau,\sigma}$ can only be nonzero if $\tau$ and $\sigma$ are on the same line, and $\sigma$ is to the right of $\tau$. The fact that there are now two lines takes into account the second part of Lemma \ref{even}, meaning the action of a central element $-\mathrm{Id}\in H_3$.
  
From this we conclude that modules $M(\mathbf{3}_{+})$, $M(\widetilde{\mathbf{3}}_{+})$ and $M(\mathbf{1}_{-})$ are also simple.

To describe $L(\mathbf{4}_{-})$, it is enough to calculate $$\mathfrak{h}^*\otimes \mathbf{4}_{-} \cong \widetilde{\mathbf{3}}_{+}\oplus \mathbf{4}_{+} \oplus \mathbf{5}_{+},$$ and  use Lemma \ref{Lemma35} to conclude 
$$L(\mathbf{4}_{-})=M(\mathbf{4}_{-})-M(\widetilde{\mathbf{3}}_{+}).$$

To describe $L(\mathbf{1}_{+})$, use Theorem \ref{support} again. The denominator of $1/5$ divides just one of the degrees of basic invariants of $H_3$, namely $10$. Thus, the support of this module is the set of all $a\in\mathfrak{h}$ whose stabilizer contains $I_2(5)$, which is a $1$-dimensional set (union of lines). That means that the character of $L(\mathbf{1}_{+})$, which is of the form 
$$\mathrm{ch}_{L(\mathbf{1}_{+})}=\mathrm{ch}_{M(\mathbf{1}_{+})}+n_{\mathbf{1}_{+},\mathbf{4}_{-}}\cdot \mathrm{ch}_{L(\mathbf{4}_{-})}+n_{\mathbf{1}_{+},\mathbf{3}_{+}}\cdot \mathrm{ch}_{L(\mathbf{3}_{+})}+n_{\mathbf{1}_{+},\widetilde{\mathbf{3}}_{+}}\cdot \mathrm{ch}_{L(\widetilde{\mathbf{3}}_{+})},$$
has a pole of order $1$ at $t=1$, i.e. that the function $$t^{-3/2}+4n_{\mathbf{1}_{+},\mathbf{4}_{-}}t^{3/2}+3n_{\mathbf{1}_{+},\mathbf{3}_{+}}t^{5/2}+3n_{\mathbf{1}_{+},\widetilde{\mathbf{3}}_{+}}t^{5/2}$$ vanishes at $t=1$ to order $2$. This translates into: $n_{\mathbf{1}_{+},\mathbf{4}_{-}}=-1$, $n_{\mathbf{1}_{+},\mathbf{3}_{+}}+ n_{\mathbf{1}_{+},\widetilde{\mathbf{3}}_{+}}=1$. This means that there is a $4$ dimensional set of singular vectors in $M(\mathbf{1}_{+})[3/2]$; using the fact that $\mathfrak{h}^*\otimes \mathbf{4}_{-} \cong \widetilde{\mathbf{3}}_{+}\oplus \mathbf{4}_{+} \oplus \mathbf{5}_{+}$, we conclude they span a copy of $L(\mathbf{4}_{-})$, so $n_{\mathbf{1}_{+},\widetilde{\mathbf{3}}_{+}}=1$, $n_{\mathbf{1}_{+},\mathbf{3}_{+}}=0$, and the Grothendieck group expression is:
$$L( \mathbf{1}_{+})=M( \mathbf{1}_{+})-M( \mathbf{4}_{-})+M( \widetilde{\mathbf{3}}_{+}).$$

To calculate the characters of $L(\mathbf{4}_{+})$, $L(\mathbf{3}_{-})$ and $L(\widetilde{\mathbf{3}}_{-})$, let us do more computations of characters of $H_3$. Namely, we use their characters to see that the multiplicity of $\mathbf{1}_{-}$ in $\mathrm{S}^3\mathfrak{h}^*\otimes \mathbf{4}_{+}$ is $1$, the multiplicity of $\mathbf{1}_{-}$ in both $\mathrm{S}^4\mathfrak{h}^*\otimes \mathbf{3}_{-}$ and $\mathrm{S}^4\mathfrak{h}^*\otimes \widetilde{\mathbf{3}}_{-}$ is $0$, and that $\mathbf{4}_{+}$ appears with multiplicity $1$ the decomposition of $\mathfrak{h}^*\otimes \widetilde{\mathbf{3}}_{-}$ and not at all in the decomposition of $\mathfrak{h}^*\otimes \mathbf{3}_{-}$. From this we can conclude: $$L(\mathbf{3}_{-})=M(\mathbf{3}_{-})$$  $$L(\mathbf{4}_{+})=M(\mathbf{4}_{+})-M(\mathbf{1}_{-})$$  $$L(\widetilde{\mathbf{3}}_{-})=M(\widetilde{\mathbf{3}}_{-})-M(\mathbf{4}_{+})+M(\mathbf{1}_{-}).$$

\section{Calculations for $c=1/3$}\label{ch1/3}

\begin{theorem}
 Irreducible representations in category $\mathcal{O}_{1/3}(H_3,\mathfrak{h})$ have the following descriptions in the Grothendieck group:
\begin{eqnarray*}
L(\mathbf{1}_{+})&=&M(\mathbf{1}_{+})-M(\mathbf{5}_{+})+M(\mathbf{4}_{-})\\
L(\mathbf{1}_{-})&=&M(\mathbf{1}_{-})\\
L(\mathbf{3}_{+})&=&M(\mathbf{3}_{+})\\
L(\mathbf{3}_{-})&=&M(\mathbf{3}_{-})\\ 
L(\widetilde{\mathbf{3}}_{+})&=&M(\widetilde{\mathbf{3}}_{+})\\
L(\widetilde{\mathbf{3}}_{-})&=&M(\widetilde{\mathbf{3}}_{-})\\ 
L(\mathbf{4}_{-})&=&M(\mathbf{4}_{-})\\
L(\mathbf{4}_{+})&=&M(\mathbf{4}_{+})-M(\mathbf{5}_{-})+M(\mathbf{1}_{-})\\
L(\mathbf{5}_{-})&=& M(\mathbf{5}_{-})-M(\mathbf{1}_{-})\\
L(\mathbf{5}_{+})&=&M(\mathbf{5}_{+})-M(\mathbf{4}_{-})
\end{eqnarray*}
None of these representations is finite dimensional.
\label{1/3}
\end{theorem}

The constants $h_{1/3}(\tau)=\frac{3}{2}-\frac{1}{3}\sum_{s\in S} s|_{\tau}$ are as follows (see Table \ref{h1/3}):
\begin{table}[h]
\begin{center}
\begin{tabular}{|c|c|c|c|c|c|c|c|c|c|} \hline
$\mathbf{1}_{+}$ & $\mathbf{1}_{-}$ & $\mathbf{3}_{+}$ & $\mathbf{3}_{-}$ & $\widetilde{\mathbf{3}}_+$ & $\widetilde{\mathbf{3}}_-$ & $\mathbf{4}_{+}$ & $\mathbf{4}_{-}$ & $\mathbf{5}_{+}$ & $\mathbf{5}_{-}$ \\ \hline
-7/2 & 13/2 & 19/6 & -1/6 & 19/6 & -1/6 & 3/2 & 3/2 & 1/2 & 5/2 \\ \hline
\end{tabular}
\end{center}
\caption{$h_{1/3}(\tau)$}
\label{h1/3}
\end{table}

Thus, $M(\mathbf{3}_{+})$, $M(\mathbf{3}_{-})$, $M(\widetilde{\mathbf{3}}_{+})$ and $M(\widetilde{\mathbf{3}}_{-})$ are simple. The other standard modules fall apart into two families:

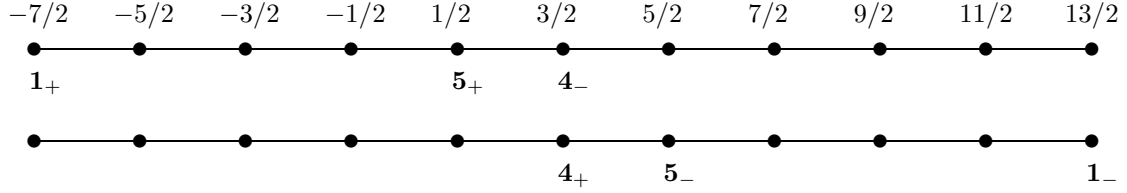
\begin{center}
\begin{picture}(400,80)
\put(0,50){\line(1,0){400}}
\put(0,15){\line(1,0){400}}

\put(0,50){\circle*{5}}
\put(40,50){\circle*{5}}
\put(80,50){\circle*{5}}
\put(120,50){\circle*{5}}
\put(160,50){\circle*{5}}
\put(200,50){\circle*{5}}
\put(240,50){\circle*{5}}
\put(280,50){\circle*{5}}
\put(320,50){\circle*{5}}
\put(360,50){\circle*{5}}
\put(400,50){\circle*{5}}
\put(0,15){\circle*{5}}
\put(40,15){\circle*{5}}
\put(80,15){\circle*{5}}
\put(120,15){\circle*{5}}
\put(160,15){\circle*{5}}
\put(200,15){\circle*{5}}
\put(240,15){\circle*{5}}
\put(280,15){\circle*{5}}
\put(320,15){\circle*{5}}
\put(360,15){\circle*{5}}
\put(400,15){\circle*{5}}
\put(-10,60){$-7/2$}
\put(30,60){$-5/2$}
\put(70,60){$-3/2$}
\put(110,60){$-1/2$}
\put(150,60){$1/2$}
\put(190,60){$3/2$}
\put(230,60){$5/2$}
\put(270,60){$7/2$}
\put(310,60){$9/2$}
\put(350,60){$11/2$}
\put(390,60){$13/2$}
\put(-2,35){$\mathbf{1}_{+}$}
\put(398,0){$\mathbf{1}_{-}$}
\put(198,35){$\mathbf{4}_{-}$}
\put(198,0){$\mathbf{4}_{+}$}
\put(158,35){$\mathbf{5}_{+}$}
\put(238,0){$\mathbf{5}_{-}$}
\end{picture}
\end{center}

So $M(\mathbf{4}_{-})$ and $M(\mathbf{1}_{-})$ are simple too.  To describe $L(\mathbf{5}_{+})$, calculate $$\mathfrak{h}\otimes \mathbf{5}_{+} \cong \mathbf{3}_{-}\oplus \widetilde{\mathbf{3}}_{-}\oplus \mathbf{4}_{-}\oplus \mathbf{5}_{-},$$
and conclude using Lemma \ref{Lemma35} $$L(\mathbf{5}_{+})=M(\mathbf{5}_{+})-M(\mathbf{4}_{-}) .$$

Next, let us describe $L(\mathbf{1}_{+})$. Again, the number of degrees of basic invariants that $3$ divides is $1$ (namely, $d_2=6$). The support is the set of all $a\in\mathfrak{h}$ that have stabilizer containing $A_2$, which is a union of lines. So, the character of $L(\mathbf{1}_{+})$ evaluated at $w=\mathrm{Id}$ has  a pole of order $1$ and the function 
$$t^{-7/2}+n_{\mathbf{1}_{+},\mathbf{5}_{+}}\cdot 5t^{1/2}+n_{\mathbf{1}_{+},\mathbf{4}_{-}}\cdot 4t^{3/2}$$ has a zero of order $2$ at $t=1$. Writing out this condition gives $n_{\mathbf{1}_{+},\mathbf{5}_{+}}=-1$, $ n_{\mathbf{1}_{+},\mathbf{4}_{-}}=1$, so the expression for $L(\mathbf{1}_{+})$ is:
$$L(\mathbf{1}_{+})=M(\mathbf{1}_{+})-M(\mathbf{5}_{+})+M(\mathbf{4}_{-}).$$

Next, we want to describe the structure of $L(\mathbf{5}_{-})$. As explained before, all the standard modules $M_{c}(\tau)$ have a contravariant bilinear form $B$ on them, whose kernel is $J_{c}(\tau)$. The form respects the grading of $M_{c}(\tau)$, in the sense that graded pieces of $M_{c}(\tau)$ are orthogonal to each other. Let the restriction of the form $B$ to $\Sym^{k}\mathfrak{h}^*\otimes \tau$ be called $B_{k}$. It is easy to compute $B_{k}$ recursively on $k$ using MAGMA algebra software. If $$L(\mathbf{5}_{-})=M(\mathbf{5}_{-})-a\cdot M(\mathbf{1}_{-}),$$ then the rank of the form $B_{4}$ on $M(\mathbf{5}_{-})$ is $$\dim L(\mathbf{5}_{-})[13/2]=\dim M(\mathbf{5}_{-})[13/2]-a\cdot \dim M(\mathbf{1}_{-})[13/2]=75-a.$$ Calculating the rank of the same $B_{4}$ in MAGMA, we get that it is $74$; hence, $a=1$ and  $$L(\mathbf{5}_{-})=M(\mathbf{5}_{-})-M(\mathbf{1}_{-}).$$

To do $L(\mathbf{4}_{+})$, notice that the multiplicity of $\mathbf{5}_{-}$ in $\Sym^{1}\mathfrak{h}^*\otimes \mathbf{4}_{+}$ is $1$, that the multiplicity of $\mathbf{1}_{-}$ in $\Sym^{5}\mathfrak{h}^*\otimes \mathbf{4}_{+}$ is $1$, and that the multiplicity of $\mathbf{1}_{-}$ in $\Sym^{4}\mathfrak{h}^*\otimes \mathbf{5}_{-}$ is $2$. So, writing out the condition that the multiplicity of $\mathbf{1}_{-}$ in $L(\mathbf{4}_{+})$ must be nonnegative, we get that the expression for it is $$L(\mathbf{4}_{+})=M(\mathbf{4}_{+})-M(\mathbf{5}_{-})+M(\mathbf{1}_{-}).$$

\section{Calculations for $c=1/2$}\label{ch1/2}

\begin{theorem}
 Irreducible representations in category $\mathcal{O}_{1/2}(H_3,\mathfrak{h})$ have the following descriptions in the Grothendieck group:
\begin{eqnarray*}
L(\mathbf{1}_{+})&=& M(\mathbf{1}_{+})-M(\mathbf{3}_{-})-M(\widetilde{\mathbf{3}}_{-})+M(\mathbf{5}_{+})-M(\mathbf{5}_{-})+M(\mathbf{3}_{+})+M(\widetilde{\mathbf{3}}_{+})-M(\mathbf{1}_{-})\\
L(\mathbf{1}_{-})&=&M(\mathbf{1}_{-})\\
L(\mathbf{3}_{+})&=&M(\mathbf{3}_{+})-M(\mathbf{1}_{-})\\
L(\mathbf{3}_{-})&=& M(\mathbf{3}_{-})-M(\mathbf{5}_{+})+M(\mathbf{5}_{-})-M(\mathbf{3}_{+})\\
L(\widetilde{\mathbf{3}}_{+})&=&M(\widetilde{\mathbf{3}}_{+})-M(\mathbf{1}_{-})\\
L(\widetilde{\mathbf{3}}_{-})&=& M(\widetilde{\mathbf{3}}_{-})-M(\mathbf{5}_{+})+M(\mathbf{5}_{-})-M(\widetilde{\mathbf{3}}_{+})\\
L(\mathbf{4}_{+})&=&M(\mathbf{4}_{+})\\
L(\mathbf{4}_{-})&=&M(\mathbf{4}_{-})\\
L(\mathbf{5}_{+})&=& M(\mathbf{5}_{+})-2\cdot M(\mathbf{5}_{-})+M(\mathbf{3}_{+})+M(\widetilde{\mathbf{3}}_{+})-M(\mathbf{1}_{-}).\\
L(\mathbf{5}_{-})&=&M(\mathbf{5}_{-})-M(\mathbf{3}_{+})-M(\widetilde{\mathbf{3}}_{+})+M(\mathbf{1}_{-})\\
\end{eqnarray*}
The following of these representations are finite dimensional: $L(\mathbf{1}_{+})$ ($\dim =115$), $L(\mathbf{3}_{-})$ (with $ch_{L(\mathbf{3}_{-})}=\chi_{\mathbf{3}_{-}}t^{-1}+\chi_{\mathbf{1}_{+}}+\chi_{\mathbf{3}_{+}}+\chi_{\mathbf{3}_{-}}t$) and $L(\widetilde{\mathbf{3}}_{-})$ (with $ch_{L(\widetilde{\mathbf{3}}_{-})}=\chi_{\widetilde{\mathbf{3}}_{-}}t^{-1}+\chi_{\mathbf{4}_{+}}+\chi_{\widetilde{\mathbf{3}}_{-}}t$).
\label{1/2}
\end{theorem}

In this case, $h_{1/2}(\tau)=\frac{3}{2}-\frac{1}{2}\sum_{s\in S} s|_{\tau}$ are (see Table \ref{h1/2}):
\begin{table}[h]
\begin{center}
\begin{tabular}{|c|c|c|c|c|c|c|c|c|c|} \hline
$\mathbf{1}_{+}$ & $\mathbf{1}_{-}$ & $\mathbf{3}_{+}$ & $\mathbf{3}_{-}$ & $\widetilde{\mathbf{3}}_+$ & $\widetilde{\mathbf{3}}_-$ & $\mathbf{4}_{+}$ & $\mathbf{4}_{-}$ & $\mathbf{5}_{+}$ & $\mathbf{5}_{-}$ \\ \hline
-6 & 9 & 4 & -1 & 4 & -1 & 3/2 & 3/2 & 0 & 3 \\ \hline
\end{tabular}
\end{center}
\caption{$h_{1/2}(\tau)$}
\label{h1/2}
\end{table}

So, $M(\mathbf{4}_{-})$ and $M(\mathbf{4}_{+})$ are simple. 

Graphic representation of Lemma \ref{even} is now
\begin{center}
\begin{picture}(300,50)
\put(0,30){\line(1,0){300}}
\put(0,30){\circle*{5}}
\put(20,30){\circle*{5}}
\put(40,30){\circle*{5}}
\put(60,30){\circle*{5}}
\put(80,30){\circle*{5}}
\put(100,30){\circle*{5}}
\put(120,30){\circle*{5}}
\put(140,30){\circle*{5}}
\put(160,30){\circle*{5}}
\put(180,30){\circle*{5}}
\put(200,30){\circle*{5}}
\put(220,30){\circle*{5}}
\put(240,30){\circle*{5}}
\put(260,30){\circle*{5}}
\put(280,30){\circle*{5}}
\put(300,30){\circle*{5}}
\put(-8,40){$-6$}
\put(92,40){$-1$}
\put(198,40){$4$}
\put(298,40){$9$}
\put(118,40){$0$}
\put(178,40){$3$}
\put(-2,15){$\mathbf{1}_{+}$}
\put(298,15){$\mathbf{1}_{-}$}
\put(98,15){$\mathbf{3}_{-}$}
\put(98,0){$\widetilde{\mathbf{3}}_{-}$}
\put(198,15){$\mathbf{3}_{+}$}
\put(198,0){$\widetilde{\mathbf{3}}_{+}$}
\put(118,15){$\mathbf{5}_{+}$}
\put(178,15){$\mathbf{5}_{-}$}
\end{picture}
\end{center}

Again, $M(\mathbf{1}_{-})$ is simple. 

Let us first analyze $L(\mathbf{3}_{-})$ and $L(\widetilde{\mathbf{3}}_{-})$. In both cases $$\dim \Hom (\mathbf{5}_{+},\mathfrak{h}^*\otimes \mathbf{3}_{-})=\dim \Hom (\mathbf{5}_{+},\mathfrak{h}^*\otimes \widetilde{\mathbf{3}}_{-})=1,$$ so by Lemma \ref{Lemma35} there is a $5$-dimensional set of lowest weight vectors at $\mathbf{h}$ weight $0$. The dimension of both of these modules at $\mathbf{h}$ weight $2$ is $$3\cdot {3+2 \choose  2}-5\cdot {2+2 \choose 2}=0.$$ This means both these modules are finite dimensional, and we can immediately determine their characters by decomposing weight spaces at $\mathbf{h}$ weights $0$ and $1$ into $H_3$ irreducible representations. They are:
$$ch_{L(\mathbf{3}_{-})}=\chi_{\mathbf{3}_{-}}t^{-1}+\chi_{\mathbf{1}_{+}}+\chi_{\mathbf{3}_{+}}+\chi_{\mathbf{3}_{-}}t$$
$$ch_{L(\widetilde{\mathbf{3}}_{-})}=\chi_{\widetilde{\mathbf{3}}_{-}}t^{-1}+\chi_{\mathbf{4}_{+}}+\chi_{\widetilde{\mathbf{3}}_{-}}t.$$

To express them in terms of characters of standard modules, write them in Grothendieck group as 
$$L(\mathbf{3}_{-})=M(\mathbf{3}_{-})-M(\mathbf{5}_{+})+n_{\mathbf{3}_{-},\mathbf{5}_{-}}M(\mathbf{5}_{-})+n_{\mathbf{3}_{-},\mathbf{3}_{+}}M(\mathbf{3}_{+})+n_{\mathbf{3}_{-},\widetilde{\mathbf{3}}_{+}}M(\widetilde{\mathbf{3}}_{-})+n_{\mathbf{3}_{-},\mathbf{1}_{-}}M(\mathbf{1}_{-})$$
$$L(\widetilde{\mathbf{3}}_{-}) = M(\widetilde{\mathbf{3}}_{-})-M(\mathbf{5}_{+})+n_{\widetilde{\mathbf{3}}_{-},\mathbf{5}_{-}}M(\mathbf{5}_{-})+n_{\widetilde{\mathbf{3}}_{-},\mathbf{3}_{+}}M(\mathbf{3}_{+})+n_{\widetilde{\mathbf{3}}_{-},\widetilde{\mathbf{3}}_{+}}M(\widetilde{\mathbf{3}}_{-})+n_{\widetilde{\mathbf{3}}_{-},\mathbf{1}_{-}}M(\mathbf{1}_{-}).$$ 
Then write the condition that dimensions of all $\mathbf{h}$ weight spaces above $2$ must be $0$ (it is enough to write the equations for weights $3$, $4$ and $9$). This produces linear equations in $n_{\tau,\sigma}$ with solutions: $n_{\mathbf{3}_{-},\mathbf{5}_{-}}=n_{\widetilde{\mathbf{3}}_{-},\mathbf{5}_{-}}=1$, $ n_{\mathbf{3}_{-},\mathbf{3}_{+}} +n_{\mathbf{3}_{-},\widetilde{\mathbf{3}}_{+}}=n_{\widetilde{\mathbf{3}}_{-},\mathbf{3}_{+}} +n_{\widetilde{\mathbf{3}}_{-},\widetilde{\mathbf{3}}_{+}}=-1, n_{\mathbf{3}_{-},\mathbf{1}_{-}}=n_{\widetilde{\mathbf{3}}_{-},\mathbf{1}_{-}}=0$. Finally, writing the $H_3$ character of $ M(\mathbf{3}_{-})[4]-M(\mathbf{5}_{+})[4]+M(\mathbf{5}_{-})[4]$ we conclude it is isomorphic to $\mathbf{3}_{+}$, whereas $ M(\widetilde{\mathbf{3}}_{-})[4]-M(\mathbf{5}_{+})[4]+M(\mathbf{5}_{-})[4]\cong \widetilde{\mathbf{3}}_{+}$, so the required expressions are: 
$$L(\mathbf{3}_{-})= M(\mathbf{3}_{-})-M(\mathbf{5}_{+})+M(\mathbf{5}_{-})-M(\mathbf{3}_{+})$$
$$L(\widetilde{\mathbf{3}}_{-})= M(\widetilde{\mathbf{3}}_{-})-M(\mathbf{5}_{+})+M(\mathbf{5}_{-})-M(\widetilde{\mathbf{3}}_{+})$$

\begin{remark}Values of $c$ for which modules with zero set of $W$ invariants exist are called \textit{aspherical}. Conjecture 4.4. in \cite{BE} was that all aspherical values are in $(-1,0)$. The module we just described, $L_{1/2}(\widetilde{\mathbf{3}}_{-})$, with character $\chi_{\widetilde{\mathbf{3}}_{-}}t^{-1}+\chi_{\mathbf{4}_{+}}+\chi_{\widetilde{\mathbf{3}}_{-}}t$, has no $H_3$ invariants and so shows that $c=1/2$ is an aspherical value for $(H_3,\mathfrak{h})$. It was the first counterexample. P.Etingof pointed out to us the notion of aspherical value and found other counterexamples in the $\mathrm{B}$ series shortly afterwards. 
\end{remark}

We use MAGMA to calculate the rank of the form $B_{5}$ on $M(\mathbf{3}_{+})$ and on $M(\widetilde{\mathbf{3}}_{+})$ and in both cases get $62$. This means there is a $3\cdot {7 \choose 2} -62=1$ dimensional kernel, and that $$L(\mathbf{3}_{+})=M(\mathbf{3}_{+})-M(\mathbf{1}_{-}),$$ $$L(\widetilde{\mathbf{3}}_{+})=M(\widetilde{\mathbf{3}}_{+})-M(\mathbf{1}_{-}).$$

To analyze $L(\mathbf{1}_{+})$, note that the number of degrees of basic invariants of $H_3$, that $2$ divides is $3$ (all the degrees $2$, $6$ and $10$ are even). This is bigger then the number of even basic invariants of any parabolic subgroup of $H_3$ except $H_3$ itself, so the support of $L(\mathbf{1}_{+})$ is the set of elements of $\mathfrak{h}$ fixed by the entire $H_3$, i.e. just a zero dimensional set consisting only of the origin. That means that 
the module $L(\mathbf{1}_{+})$ is finite dimensional.

\begin{remark}
Notice that the previous argument depended only on the denominator of $c=1/2$; it actually proves that $L_{r/2}(\mathbf{1}_{+})$ is finite dimensional for all odd $r\ge 0$.
\end{remark}

Now we use MAGMA \cite{BCP} to calculate the rank of the form $B$ restricted to $M(\mathbf{1}_{+})[-1]=\mathrm{S}^5\mathfrak{h}^*\otimes \mathbf{1}_{+}$. This is  $21$ dimensional space, and the rank of the form is $15$. Since both $\mathbf{3}_{-}$ and $\widetilde{\mathbf{3}}_{-}$ appear in the decomposition of $\mathrm{S}^5\mathfrak{h}^*$ into $H_3$ subrepresentations, and each of them with multiplicity $2$, we need some more calculations to see how this $6$ dimensional space of singular vectors looks. To do that, again use MAGMA to compute the $H_3$ character on the $6$ dimensional kernel of $B$ on $M(\mathbf{1}_{+})[-1]$. This computation shows that the kernel is $\mathbf{3}_{-}\oplus \widetilde{\mathbf{3}}_{-}$, so $n_{\mathbf{1}_{+}, \mathbf{3}_{-}}=n_{\mathbf{1}_{+}, \widetilde{\mathbf{3}}_{-}}=1$. 

Because $L(\mathbf{1}_{+})$ is finite dimensional and has an $\mathfrak{sl}_2$ representation structure, we know that $\dim L(\mathbf{1}_{+})[j]=\dim L(\mathbf{1}_{+})[-j]$ for every integer $j$. This gives us a system of linear equations whose only solution yields the following expression for the irreducible module:
$$L(\mathbf{1}_{+})=M(\mathbf{1}_{+})-M(\mathbf{3}_{-})-M(\widetilde{\mathbf{3}}_{-})+M(\mathbf{5}_{+}) -M(\mathbf{5}_{-})+n_{\mathbf{1}_{+},\mathbf{3}_{+}}M(\mathbf{3}_{+})+n_{\mathbf{1}_{+},\widetilde{\mathbf{3}}_{+}}M(\widetilde{\mathbf{3}}_{+})-M(\mathbf{1}_{-})$$
with $n_{\mathbf{1}_{+},\mathbf{3}_{+}}+n_{\mathbf{1}_{+},\widetilde{\mathbf{3}}_{+}}=2.$

To calculate $n_{\mathbf{1}_{+},\mathbf{3}_{+}},n_{\mathbf{1}_{+},\widetilde{\mathbf{3}}_{+}}$ we make the following observation. As the copy of $\mathfrak{sl}_{2}$ in $H_{c}(H_3,\mathfrak{h})$ commutes with $H_3$, for any $H_{c}(H_3,\mathfrak{h})$ module $M$ and any irreducible representation $\tau$ of $H_3$ we can put the $\mathfrak{sl}_{2}$ module structure on $\Hom_{H_3}(\tau,M)$ by letting  $\mathfrak{sl}_{2}$ act on the value. If $M=L(\mathbf{1}_{+})$, then this module is finite dimensional, so dimensions of weight spaces are symmetric around $0$. In other words, $\dim \Hom_{H_3}(\tau, L(\mathbf{1}_{+})[-j])=\dim \Hom_{H_3}(\tau, L(\mathbf{1}_{+})[j])$. Doing this computation for $\tau=\mathbf{3}_{+}$ and $j=4$ gives us that this dimension is $\dim \Hom_{H_3}(\mathbf{3}_{+},\Sym^2\mathfrak{h}^*)=0$. 

Representations of $A_{5}$ and $H_3$ are defined over the field $\mathbb{Q}[\sqrt{5}]$, which is a field extension of $\mathbb{Q}$ of degree $2$. The Galois action of $\mathbb{Z}_{2}$ corresponding to this extension is $\sqrt{5}\mapsto -\sqrt{5}$. It acts on all characters, and it is clear from the character table \ref{H3CharTable} that the action of the Galois group on the character of a representation $V$ of $H_{3}$ is trivial if and only if $$\dim \Hom_{H_{3}}(\mathbf{3}_{-}, V)+ \dim \Hom_{H_{3}}(\mathbf{3}_{+}, V)=\dim \Hom_{H_{3}}(\widetilde{\mathbf{3}}_{-}, V)+\dim \Hom_{H_{3}}(\widetilde{\mathbf{3}}_{+}, V),$$ in other words, if, seen as a representation of $A_{5}$ and decomposed into irreducible subrepresentations, $V$ has the same multiplicity of $\mathbf{3}$ and $\widetilde{\mathbf{3}}$.

Calculation of the $H_3$ characters for $L(\mathbf{1}_{+})[-4]$ and $L(\mathbf{1}_{+})[4]=\Sym^{10}\mathfrak{h}^*\otimes \mathbf{1}_{+}-\Sym^{5}\mathfrak{h}^*\otimes (\mathbf{3}_{+}\oplus \widetilde{\mathbf{3}}_{+})+\Sym^{4}\mathfrak{h}^*\otimes \mathbf{5}_{+}-\Sym^{1}\mathfrak{h}^*\otimes \mathbf{5}_{-}+n_{\mathbf{1}_{+},\mathbf{3}_{+}}\mathbf{3}_{+}+n_{\mathbf{1}_{+}, \widetilde{\mathbf{3}}_{+}}\widetilde{\mathbf{3}}_{+}$ (an elementary computation of $H_3$ characters, though a tedious one) show the character of $L(\mathbf{1}_{+})[-4]$ is invariant under the above Galois action, and that the character of  $L(\mathbf{1}_{+})[4]$ (which is the same) is invariant if and only if $n_{\mathbf{1}_{+},\mathbf{3}_{+}}=n_{\mathbf{1}_{+}, \widetilde{\mathbf{3}}_{+}}$. So, they both have to be $1$, and the character is, as claimed in the theorem,
$$L(\mathbf{1}_{+})=M(\mathbf{1}_{+})-M(\mathbf{3}_{-})-M(\widetilde{\mathbf{3}}_{-})+M(\mathbf{5}_{+}) -M(\mathbf{5}_{-})+M(\mathbf{3}_{+})+M(\widetilde{\mathbf{3}}_{+})-M(\mathbf{1}_{-}).$$

\subsection{Cherednik algebra $H_{1/2}(\mathbb{Z}_{2}\times\mathbb{Z}_{2},\mathfrak{h}')$ and calculation of $L(\mathbf{5}_{-})$}
We will calculate $L(\mathbf{5}_{-})$ using the induction functor, see \ref{indchap}. To do that, let us first describe the algebra we will be inducing from.

A way to get a maximal parabolic subgroups of Coxeter groups is to remove one vertex from the Coxeter graph, which corresponds to removing one generator. In this case, let us remove the middle vertex of the $H_3$ graph, and thus get a disconnected graph 
\begin{center}
\begin{picture}(60,30)
\put(0,10){\circle*{5}}
\put(60,10){\circle*{5}}
\end{picture}
\end{center}
 of $\mathbb{Z}_{2}\times \mathbb{Z}_{2}$. In the isomorphism $H_3\cong \mathbb{Z}_{2}\times A_{5}$, we can take the Coxeter generators of $H_{3}$ to be $s_{1}=-(12)(34), s_{2}=-(15)(34), s_3=-(13)(24)$. Then the generators of $W'=\mathbb{Z}_{2}\times \mathbb{Z}_{2}$ are $s_{1},s_{3}$. Let us write the character table of $\mathbb{Z}_{2}\times \mathbb{Z}_{2}$, with the main purpose of introducing notation and names of representations:  see Table \ref{Z2xZ2 character}.

\begin{table}[h!]
\begin{center}
\begin{tabular}{|c||c|c|c|c|} \hline
& $\mathrm{Id}$ & $-(12)(34)$ & $-(13)(24)$ & $(14)(23)$ \\ \hline
$\mathbf{1}_{++}$ & 1 & 1 & 1 & 1 \\ \hline
$\mathbf{1}_{+-}$ & 1 & 1 & -1 & -1 \\ \hline
$\mathbf{1}_{-+}$ & 1 & -1 & 1 & -1 \\ \hline
$\mathbf{1}_{--}$ & 1 & -1 & -1 & 1 \\ \hline
\end{tabular}
\end{center}
\caption{Character table for $\Bbb{Z}_2 \times \Bbb{Z}_2$}
\label{Z2xZ2 character}
\end{table} 

Working out the irreducible modules $L_{1/2}(\tau), \tau \in \hat{W}'$ is really easy in this case. They have the lowest weights given in table \ref{h1/2,z2}.
\begin{table}[h!]
\begin{center}
\begin{tabular}{|c|c|c|c|} \hline
$\mathbf{1}_{++}$ & $\mathbf{1}_{+-}$ & $\mathbf{1}_{-+}$ & $\mathbf{1}_{--}$ \\ \hline
0 & 1 & 1 & 2  \\ \hline
\end{tabular}
\end{center}
\caption{$h_{1/2}(\tau),  \tau$ irreducible representation of $\mathbb{Z}_2\times\mathbb{Z}_2$}
\label{h1/2,z2}
\end{table}
So using only Lemma \ref{Lemma35} and the fact $\mathfrak{h}'\cong \mathbf{1}_{-+}\oplus \mathbf{1}_{+-}$, we get that the module $$L_{1/2}(\mathbf{1}_{++})=M_{1/2}(\mathbf{1}_{++})-M_{1/2}(\mathbf{1}_{+-})-M_{1/2}(\mathbf{1}_{-+})+M_{1/2}(\mathbf{1}_{--})$$ is a one dimensional representation of $H_{1/2}(\mathbb{Z}_{2}\times\mathbb{Z}_{2},\mathfrak{h}')$.

Let $b$ be any point whose stabilizer is this copy of $\mathbb{Z}_{2}\times\mathbb{Z}_{2}$. We are going to apply the induction functor $\Ind_{b}$ to the one dimensional module $L_{1/2}(\mathbf{1}_{++})$. Before we do that, let us decompose all the representations of $H_{3}$ into representations of $\mathbb{Z}_{2}\times \mathbb{Z}_{2}$:
\begin{table}[h!]
\begin{center}
\begin{tabular}{|c||c|c|c|c|c|} \hline
& $\mathrm{Id}$ & $-(12)(34)$ & $-(13)(24)$ & $(14)(23)$ & $\cong$ \\ \hline
$\mathbf{1}_{+}$ & 1 & 1 & 1 & 1 & $\mathbf{1}_{++}$ \\ \hline
$\mathbf{1}_{-}$ & 1 & -1 & -1 & 1 & $\mathbf{1}_{--}$ \\ \hline
$\mathbf{3}_{+}$ & 3 & -1 & -1 & -1  & $\mathbf{1}_{-+}\oplus \mathbf{1}_{+-}\oplus \mathbf{1}_{--}$  \\ \hline
$\mathbf{3}_{-}$ & 3 & 1 & 1 & -1  &  $\mathbf{1}_{-+}\oplus \mathbf{1}_{+-}\oplus \mathbf{1}_{++}$ \\ \hline
$\widetilde{\mathbf{3}}_{-}$ & 3 & -1 & -1 & -1 &  $\mathbf{1}_{-+}\oplus \mathbf{1}_{+-}\oplus \mathbf{1}_{--}$ \\ \hline
$\widetilde{\mathbf{3}}_{-}$ & 1 & 1 & 1 & -1 &  $\mathbf{1}_{-+}\oplus \mathbf{1}_{+-}\oplus \mathbf{1}_{++}$ \\ \hline
$\mathbf{4}_{+}$ & 4 & 0 & 0 & 0  &  $\mathbf{1}_{--}\oplus \mathbf{1}_{-+}\oplus \mathbf{1}_{+-}\oplus \mathbf{1}_{++}$ \\ \hline
$\mathbf{4}_{-}$ & 4 & 0 & 0 & 0  & $\mathbf{1}_{--}\oplus \mathbf{1}_{-+}\oplus \mathbf{1}_{+-}\oplus \mathbf{1}_{++}$ \\ \hline
$\mathbf{5}_{+}$ & 5 & 1 & 1 & 1  &$\mathbf{1}_{--}\oplus \mathbf{1}_{-+}\oplus \mathbf{1}_{+-}\oplus 2\cdot \mathbf{1}_{++}$ \\ \hline
$\mathbf{5}_{-}$ & 5 & -1 & -1 & -1  & $2\cdot \mathbf{1}_{--}\oplus \mathbf{1}_{-+}\oplus \mathbf{1}_{+-}\oplus \mathbf{1}_{++}$ \\ \hline
\end{tabular}
\end{center}
\caption{Decomposition of irreducible representations of $H_3$ as representations of $\mathbb{Z}_{2}\times\mathbb{Z}_2$}
\label{Z2xZ2 subrep}
\end{table} 

So, using Proposition \ref{ind}, the expression in the Grothendieck group of $\mathcal{O}_{1/2}(H_{3},\mathfrak{h})$ for the induced module $\Ind _{b}(L_{1/2}(\mathbf{1}_{++}))$ is 
\begin{eqnarray*}\Ind _{b}(L_{1/2}(\mathbf{1}_{++}))&=&\Ind _{b}(M_{1/2}(\mathbf{1}_{++}))-\Ind _{b}(M_{1/2}(\mathbf{1}_{-+})) -\Ind _{b}(M_{1/2}(\mathbf{1}_{+-}))+\Ind _{b}(L_{1/2}(\mathbf{1}_{--}))\\
&=& M(\mathbf{1}_{+})-M(\mathbf{3}_{+})-M(\widetilde{\mathbf{3}}_{+})+M(\mathbf{5}_{-})+M(\mathbf{5}_{+})-M(\mathbf{3}_{-})-M(\widetilde{\mathbf{3}}_{-})+M(\mathbf{1}_{-}).\\
\end{eqnarray*}

This means there is a module in $\mathcal{O}_{1/2}(H_{3},\mathfrak{h})$ with this expression in the Grothendieck group. Its composition series must contain an irreducible module containing $M(\mathbf{1}_{+})$ in its Grothendieck group expression, and there is only one such. Subtracting the known Grothendieck group expression of $L(\mathbf{1}_{+})$ from the one for $\Ind _{b}(L_{1/2}(\mathbf{1}_{++}))$, we get that there must exist a module with Grothendieck group expression $$2\big( M(\mathbf{5}_{-})-M(\mathbf{3}_{+})-M(\widetilde{\mathbf{3}}_{+})+M(\mathbf{1}_{-})\big) .$$

Now, Lemma \ref{Lemma35} and the decomposition of $\mathfrak{h}^*\otimes\mathbf{5}_{-}$ into irreducible subrepresentations imply that the irreducible module $L(\mathbf{5}_{-})$ is of the form $$L(\mathbf{5}_{-})=M(\mathbf{5}_{-})-M(\mathbf{3}_{+})-M(\widetilde{\mathbf{3}}_{+})+a\cdot M(\mathbf{1}_{-}),$$ with $a\in \mathbb{Z}$.  There are $3$ copies of $\mathbf{1}_{-}$ in $M(\mathbf{5}_{-})[9]$, 2 copies of $\mathbf{1}_{-}$ in $M(\mathbf{3}_{+})[9]$ and 2 copies of $\mathbf{1}_{-}$ in $M(\widetilde{\mathbf{3}}_{+})[9]$, so $3-2-2+a\ge 0$ and $a\ge 1$.

Subtracting two times this expression from the above expression for the module we concluded must exist, we get that there also must be a module with Grothendieck group expression
$$2(1-a)M(\mathbf{1}_{-}),$$ i.e. that $a\le1$, so $a=1$. This proves that the expression for the irreducible module we wanted is $$L(\mathbf{5}_{-})=M(\mathbf{5}_{-})-M(\mathbf{3}_{+})-M(\widetilde{\mathbf{3}}_{+})+M(\mathbf{1}_{-}).$$

\subsection{Cherednik algebra $H_{1/2}(S_{3},\mathfrak{h}')$ and calculation of $L(\mathbf{5}_{+})$}

We start by doing the MAGMA computation of rank of $B$ in degrees $3$ and $4$ we get that it is $40$ and $51$, so the Grothendieck group expression is of the form
$$L(\mathbf{5}_{+})=M(\mathbf{5}_{+})-2\cdot M(\mathbf{5}_{-})+n_{\mathbf{5}_{-},\mathbf{3}_{+}}M(\mathbf{3}_{+})+n_{\mathbf{5}_{-},\widetilde{\mathbf{3}}_{+}}M(\widetilde{\mathbf{3}}_{+})+n_{\mathbf{5}_{-},\mathbf{1}_{-}}M(\mathbf{1}_{-}),$$ with $n_{\mathbf{5}_{-},\widetilde{\mathbf{3}}_{+}}+n_{\mathbf{5}_{-},\mathbf{3}_{+}}=2$. Looking at the dimension of $L(\mathbf{5}_{+})[k]$, which is a quadratic polynomial in $k$ with leading term $\frac{1}{2}(1+n_{\mathbf{5}_{-},\mathbf{1}_{-}})k^2$, and writing the condition that it is $\ge 0$ for large $k$, we conclude $n_{\mathbf{5}_{-},\mathbf{1}_{-}} \ge -1$.

To finish the analysis, we need to look at another Cherednik algebra associated to a parabolic subgroup of $H_3$, like in the last section. This time, remove the rightmost vertex in the Coxeter graph, to get a group $W'=S_3$ generated by $s_{1}=-(12)(34)$, $s_{2}=-(15)(34)$. Its character table is very well known:

\begin{table}[h!]
\begin{center}
\begin{tabular}{|c||c|c|c|} \hline
& $\mathrm{Id}$ & $s_{1}$ & $(125)$\\ \hline
\# & $1$ & $3$ & $2$\\ \hline
$\mathbf{1}_{+}$ & 1 & 1 & 1  \\ \hline
$\mathbf{1}_{-}$ & 1 & -1 & 1  \\ \hline
$\mathbf{2}$ & 2 & 0 & -1  \\ \hline
\end{tabular}
\end{center}
\caption{Character table for $S_3$}
\label{S3 character}
\end{table} 

Working out the irreducible modules $L_{1/2}(\tau), \tau \in \hat{W}'$ is again really easy. The lowest weights are
\begin{table}[h]
\begin{center}
\begin{tabular}{|c|c|c|} \hline
 $\mathbf{1}_{+}$ & $\mathbf{1}_{-}$ & $\mathbf{2}$ \\ \hline
-1/2 & 5/2 & 1  \\ \hline
\end{tabular}
\end{center}
\caption{$h_{1/2}(\tau), \tau$ irreducible representation of $S_3$}
\label{h1/2,s3}
\end{table}

Because the denominator of $1/2$ is a degree of a basic invariant of $S_{3}$, the category $\mathcal{O}_{1/2}(S_{3},\mathfrak{h}')$ cannot be semisimple. So, $M_{1/2}(\mathbf{1}_{+})$ is not simple (as the other two are). Looking at the possible options and decomposing $\Sym^{2} \mathfrak{h}'\otimes \mathbf{1}_{+}=\Sym^{2}\mathbf{2}=\mathbf{2}\oplus \mathbf{1}_{-}$, we conclude $L_{1/2}(\mathbf{1}_{+})=M_{1/2}(\mathbf{1}_{+})-M_{1/2}(\mathbf{1}_{-})$.

Now let $b$ be any point with a stabilizer $W'$ and apply $\Ind_{b}$ to $L_{1/2}(\mathbf{1}_{+})$. In the same way as before (using decompositions of $H_3$ representations into $S_3$ irreducible components, and applying Lemma \ref{ind}), we get that there is a module in $\mathcal{O}_{c}(H_3,\mathfrak{h})$ with the Grothendieck group description
$$M(\mathbf{1}_{+})+M(\mathbf{3}_{-})+M(\widetilde{\mathbf{3}}_{-})+M(\mathbf{5}_{+})-M(\mathbf{5}_{-})-M(\mathbf{3}_{+})-M(\widetilde{\mathbf{3}}_{+})-M(\mathbf{1}_{-}).$$
Subtracting $L(\mathbf{1}_{+})$, which has to be in its composition series, from it, and doing the same thing for $L(\mathbf{3}_{-})$ and  $L(\widetilde{\mathbf{3}}_{-})$, we see that there is a module in $\mathcal{O}_{c}(H_3,\mathfrak{h})$ with the Grothendieck group description $$4(M(\mathbf{5}_{+})-M(\mathbf{5}_{-})).$$

$L(\mathbf{5}_{+})$ must appear as a factor in the composition series of this module $4$ times. So, subtract $4\cdot L(\mathbf{5}_{+})$ from it to get that there exists a module with  expression $$4\cdot M(\mathbf{5}_{-})-4 \cdot n_{\mathbf{5}_{-},\mathbf{3}_{+}}M(\mathbf{3}_{+})-4 \cdot n_{\mathbf{5}_{-},\widetilde{\mathbf{3}}_{+}}M(\widetilde{\mathbf{3}}_{+})-4 \cdot n_{\mathbf{5}_{-},\mathbf{1}_{-}}M(\mathbf{1}_{-}).$$ Subtracting the known expression for $4 \cdot L(\mathbf{5}_{-})$, we get that there must be a module with expression: $$4 \cdot(1- n_{\mathbf{5}_{-},\mathbf{3}_{+}}) M(\mathbf{3}_{+})+4\cdot  (1-n_{\mathbf{5}_{-},\widetilde{\mathbf{3}}_{+}})M(\widetilde{\mathbf{3}}_{+})-4\cdot  (1+ n_{\mathbf{5}_{-},\mathbf{1}_{-}})M(\mathbf{1}_{-}).$$

This implies $1- n_{\mathbf{5}_{-},\mathbf{3}_{+}}\ge 0$ and $1- n_{\mathbf{5}_{-},\widetilde{\mathbf{3}}_{+}}\ge 0$, which together with $n_{\mathbf{5}_{-},\mathbf{3}_{+}}+n_{\mathbf{5}_{-},\widetilde{\mathbf{3}}_{+}}=2$ means $n_{\mathbf{5}_{-},\mathbf{3}_{+}}=n_{\mathbf{5}_{-},\widetilde{\mathbf{3}}_{+}}=1$. The last module then becomes $$-4\cdot  (1+ n_{\mathbf{5}_{-},\mathbf{1}_{-}})M(\mathbf{1}_{-}), $$ so $1+ n_{\mathbf{5}_{-},\mathbf{1}_{-}}\le 0$, which means $n_{\mathbf{5}_{-},\mathbf{1}_{-}}= -1$. Therefore, we have $$L(\mathbf{5}_{+})=M(\mathbf{5}_{+})-2\cdot M(\mathbf{5}_{-})+M(\mathbf{3}_{+})+M(\widetilde{\mathbf{3}}_{+})-M(\mathbf{1}_{-}).$$

\begin{remark}
As explained before, for $c=r/d$, $d\ge 3$, we will use equivalences of categories $\mathcal{O}_{c}\to \mathcal{O}_{rc}$ to get the descriptions of all modules $L_{rc}(\tau)$. These functors are not available in the case of $d=2$, so here we use different equivalences $S_{r/2}^{+}:\mathcal{O}_{r/2}\to \mathcal{O}_{(r+2)/2}$. However, $S_{c}^{+}$ is only an equivalence of categories when $\mathcal{O}_{c}^{+}=0$. This is not the case for $c=1/2$, as we have seen an example of a module $L_{1/2}(\widetilde{\mathbf{3}}_{-})$ that has a zero set of $H_{3}$ invariants, therefore being in $\mathcal{O}_{c}^{+}$, with $S_{1/2}^+(L_{1/2}(\widetilde{\mathbf{3}}_{-}))=0$. This is why we cannot use results for $1/2$ to derive results for $r/2$ for any positive odd $r$.
\end{remark}

\section{Calculations for $c=3/2$}\label{ch3/2}

\begin{theorem}
 Irreducible representations in category $\mathcal{O}_{3/2}(H_3,\mathfrak{h})$ have the following descriptions in the Grothendieck group:
\begin{eqnarray*}
L(\mathbf{1}_{+})&=& M(\mathbf{1}_{+})-M(\mathbf{3}_{-})-M(\widetilde{\mathbf{3}}_{-})+M(\mathbf{5}_{+})- M(\mathbf{5}_{-})+M(\mathbf{3}_{+})+M(\widetilde{\mathbf{3}}_{+})-M(\mathbf{1}_{-})\\
L(\mathbf{1}_{-})&=&M(\mathbf{1}_{-})\\
L(\mathbf{3}_{+})&=&M(\mathbf{3}_{+})-M(\mathbf{1}_{-})\\
L(\mathbf{3}_{-})&=& M(\mathbf{3}_{-})-M(\mathbf{5}_{+})+M(\mathbf{5}_{-})-M(\mathbf{3}_{+})\\
L(\widetilde{\mathbf{3}}_{+})&=&M_{1/2}(\widetilde{\mathbf{3}}_{+})-M_{1/2}(\mathbf{1}_{-})\\
L(\widetilde{\mathbf{3}}_{-})&=& M(\widetilde{\mathbf{3}}_{-})-M(\mathbf{5}_{+})+M(\mathbf{5}_{-})-M(\widetilde{\mathbf{3}}_{+})\\
L(\mathbf{4}_{+})&=&M(\mathbf{4}_{+})\\
L(\mathbf{4}_{-})&=&M(\mathbf{4}_{-})\\
L(\mathbf{5}_{+})&=& M(\mathbf{5}_{+})-2\cdot M(\mathbf{5}_{-})+M(\mathbf{3}_{+})+M(\widetilde{\mathbf{3}}_{+})-M(\mathbf{1}_{-})\\
L(\mathbf{5}_{-})&=&M(\mathbf{5}_{-})-M(\mathbf{3}_{+})-M(\widetilde{\mathbf{3}}_{+})+M(\mathbf{1}_{-})\\
\end{eqnarray*}
Three of these representations are finite dimensional: $L(\mathbf{1}_{+})$, $L(\mathbf{3}_{-})$ and $L(\widetilde{\mathbf{3}}_{-})$. 
\label{3/2}
\end{theorem}

In this case, $h_{3/2}(\tau)=\frac{3}{2}-\frac{1}{2}\sum_{s\in S} s|_{\tau}$ are (see Table \ref{h1/2}):
\begin{table}[h]
\begin{center}
\begin{tabular}{|c|c|c|c|c|c|c|c|c|c|} \hline
$\mathbf{1}_{+}$ & $\mathbf{1}_{-}$ & $\mathbf{3}_{+}$ & $\mathbf{3}_{-}$ & $\widetilde{\mathbf{3}}_+$ & $\widetilde{\mathbf{3}}_-$ & $\mathbf{4}_{+}$ & $\mathbf{4}_{-}$ & $\mathbf{5}_{+}$ & $\mathbf{5}_{-}$ \\ \hline
-21 & 24 & 9 & -6 & 9 & -6 & 3/2 & 3/2 & -3 & 6 \\ \hline
\end{tabular}
\end{center}
\caption{$h_{3/2}(\tau)$}
\label{h3/2}
\end{table}

Graphic representation of Lemma \ref{even} is:
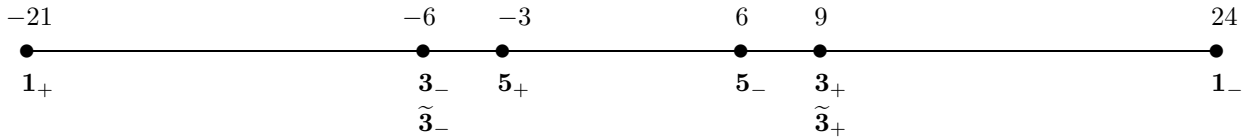
\begin{center}
\begin{picture}(450,50)
\put(0,30){\line(1,0){450}}

\put(0,30){\circle*{5}}
\put(150,30){\circle*{5}}
\put(180,30){\circle*{5}}
\put(270,30){\circle*{5}}
\put(300,30){\circle*{5}}
\put(450,30){\circle*{5}}
\put(-8,40){$-21$}
\put(142,40){$-6$}
\put(178,40){$-3$}
\put(268,40){$6$}
\put(298,40){$9$}
\put(448,40){$24$}
\put(-2,15){$\mathbf{1}_{+}$}
\put(148,15){$\mathbf{3}_{-}$}
\put(148,0){$\widetilde{\mathbf{3}}_{-}$}
\put(178,15){$\mathbf{5}_{+}$}
\put(268,15){$\mathbf{5}_{-}$}
\put(298,15){$\mathbf{3}_{+}$}
\put(298,0){$\widetilde{\mathbf{3}}_{+}$}
\put(448,15){$\mathbf{1}_{-}$}
\end{picture}
\end{center}

First, because $S_{1/2}^+$ is an equivalence of $\mathcal{O}_{1/2}/\mathcal{O}_{1/2}^+$ and $\mathcal{O}_{3/2}/\mathcal{O}_{3/2}^-$, we can conclude that the modules $L(\tau)=L_{3/2}(\tau)$ for $\tau \in \{ \mathbf{4}_{+}, \mathbf{4}_{-}, \mathbf{5}_{+}, \mathbf{5}_{-}, \mathbf{3}_{+}, \widetilde{\mathbf{3}}_{+},\mathbf{1}_{-}\}$ have the Grothendieck group expressions analogous to those in $c=1/2$ case.

The argument from the previous chapter shows that $L(\mathbf{1}_{+})$ is finite dimensional. Calculating the rank of the form $B_{15}$ and the character on the kernel lets us conclude that $J(\mathbf{1}_{+})[-6]\cong \mathbf{3}_{-}\oplus \widetilde{\mathbf{3}_{-}}$. Solving the system of equations $\dim L(\mathbf{1}_{+})[k]=\dim L(\mathbf{1}_{+})[-k]$ in $n_{\mathbf{1}_{+},\sigma}$ (it is enough to do so for $k=3,6,9,24$) gives all the coefficients of the Grothendieck group, except $n_{\mathbf{1}_{+},\mathbf{3}_{+}}$ and $n_{\mathbf{1}_{+},\widetilde{\mathbf{3}}_{+}}$, for which we can only conclude that their sum is $2$. Then we look at the $H_3$ characters on spaces $L(\mathbf{1}_{+})[9]$ and $L(\mathbf{1}_{+})[-9]$; the condition that they must be invariant under the Galois group action $\sqrt{5}\mapsto -\sqrt{5}$ implies that $n_{\mathbf{1}_{+},\mathbf{3}_{+}}=n_{\mathbf{1}_{+},\widetilde{\mathbf{3}}_{+}}$. This gives us the desired formula for $L(\mathbf{1}_{+})$ in terms of $M(\sigma)$.

To describe $L(\mathbf{3}_{-})$ and $L(\widetilde{\mathbf{3}}_{+})$, we first use MAGMA to compute the rank of the form $B$ on the $\mathbf{h}$ weight space $-3$. In both cases we get that there is a $5$ dimensional kernel. Writing out the dimensions of graded pieces we can again conclude that these modules are finite dimensional, with the Grothendieck group expressions
$$L(\mathbf{3}_{-})=M(\mathbf{3}_{-})-M(\mathbf{5}_{+})+M(\mathbf{5}_{-})+n_{\mathbf{3}_{-},\mathbf{3}_{+}}M(\mathbf{3}_{+})+n_{\mathbf{3}_{-},\widetilde{\mathbf{3}}_{+}}M(\widetilde{\mathbf{3}}_{+})$$
$$L(\widetilde{\mathbf{3}}_{-})=M(\widetilde{\mathbf{3}}_{-})-M(\mathbf{5}_{+})+M(\mathbf{5}_{-})+n_{\widetilde{\mathbf{3}}_{-},\mathbf{3}_{+}}M(\mathbf{3}_{+})+n_{\widetilde{\mathbf{3}}_{-},\widetilde{\mathbf{3}}_{+}}M(\widetilde{\mathbf{3}}_{+}),$$ with $n_{\mathbf{3}_{-},\mathbf{3}_{+}}+n_{\mathbf{3}_{-},\widetilde{\mathbf{3}}_{+}}=n_{\widetilde{\mathbf{3}}_{-},\mathbf{3}_{+}}+n_{\widetilde{\mathbf{3}}_{-},\widetilde{\mathbf{3}}_{+}}=-1$. Finally, looking at the trace of an element $(12345)$ on $L(\mathbf{3}_{-})$ and $L(\widetilde{\mathbf{3}}_{-})$, which of course needs to be $0$, we can conclude that all $n_{\tau, \sigma}$ are as in the statement of the theorem.

The representation $L(\widetilde{\mathbf{3}}_{-})$ contains no $H_{3}$ anti-invariants, i.e. no copy of $\mathbf{1}_{-}$. This is the representation that $S_{1/2}^{-}$ annihilates. While the functors $S_{1/2}^+$ and $S_{1/2}^{-}$ are not equivalences of categories, they are if we work modulo these two irreducible representations. However, the module  $L(\widetilde{\mathbf{3}}_{-})$ has an $H_3$ invariant, for example a one dimensional space in degree $3$. All other modules also have nontrivial $H_3$ invariants. So we can conclude that $S_{3/2}^+$ is an equivalence of categories, and that all the irreducible modules in $\mathcal{O}_{5/2}$ will have the Grothendieck group expressions equivalent to corresponding modules in $\mathcal{O}_{3/2}$. It is easy to see that then all these modules are going to have a nonzero $H_3$ invariant (because all the modules for $c=3/2$ have them in some low degree, and $L_{5/2}(\tau)[k]=M_{5/2}(\tau)[k]$ for at least as many degrees $k$ as was the case for $c=3/2$). So, we can conclude:

\begin{lemma} Functors $S^{+}_{c}:\mathcal{O}_{c}\to \mathcal{O}_{c+1}$ are equivalences of categories for $c=r/2$, $r$ odd, $r\ge 3$. \label{equi3}
\end{lemma}

This lemma allows us to derive formulas for Grothendieck group expressions of $L_{c}(\tau)$ in terms of $M_{c}(\tau)$ for all $c=r/2, r>3$. It is used in the proof of Theorem \ref{main}.


\begin{thebibliography}{99}
\begin{normalsize}

\bibitem[BCP]{BCP} Bosma, W, Cannon J, and Playoust, C, \emph{The Magma algebra system. I. The user language}, J. Symb. Comput. Vol. 24, No. 3-4, (1997) 235-265. http://magma.maths.usyd.edu.au/magma/

\bibitem[BE]{BE} R. Bezrukavnikov and P. Etingof, {\it Parabolic induction and restriction functors for rational Cherednik algebras}, Sel. Math. Vol. 14, No. 3-4 (2009) 397-425.

\bibitem[BEG1]{BEG1}  Y. Berest , P. Etingof and V. Ginzburg, {\it Finite dimensional representations of rational Cherednik algebras}, Int. Math. Res. Not. 19 (2003), 1053-1088.

\bibitem[BEG2]{BEG2}  Y. Berest , P. Etingof and V. Ginzburg, {\it Cherednik algebras and differential operators on quasi-invariants},  Duke Math. J. Vol. 118, No 2 (2003), 279-337. 

\bibitem[Ch]{Ch} I. Cherednik, {\it Double affine Hecke algebras and Macdonald conjectures}, Ann. of Math. (2) Vol 141, No 1. (1995), 191-216. 

\bibitem[Chm]{Chm}  T. Chmutova, {\it Representations of the rational Cherednik algebras of dihedral type}, J. of Algebra Vol. 297, No. 2 (2006), 542-565.

\bibitem[D]{D} C.F. Dunkl, {\it Differential-difference operators associated to reflection groups}, Trans. Amer. Math. Soc, Vol. 311, No. 1 (1989), 167-183 .

\bibitem[DO]{DO} C.F. Dunkl and E.M. Opdam, {\it Dunkl opeartors for complex reflection groups}, Proc. of the Lond. Math. Soc. Vol. 86 (2003), 70-108.

\bibitem[DJO]{DJO} C.F. Dunkl, M.F.E. de Jeu and E.M. Opdam, {\it Singular polynomials for finite reflection groups}, Trans. Amer. Math. Soc. Vol. 346 (1994), 237-256.

\bibitem[E]{E} P. Etingof, {\it Lectures on Calogero-Moser systems and representation theory}, Zurich lectures in advanced mathematics, European mathematical society, Z\"urcih 2007 

\bibitem[E2]{E2} P. Etingof, {\it Supports of irreducible spherical representations of rational Cherednik algebras of finite Coxeter groups}, preprint arXiv:0911.3208

\bibitem[EG]{EG} P. Etingof and V. Ginzburg, {\it Symplectic reflection algebras, Calogero-Moser space, and deformed Harish-Chandra isomorphism}, Invent. Math. Vol. 147, No. 2, (2002), 243-348.

\bibitem[EM]{EM} P. Etingof, X. Ma, {\it Lecture notes on Cherednik algebras}, math.mit.edu/~etingof

\bibitem[ES]{ES} P. Etingof  and E. Stoica, {\it Unitary representations of rational Cherednik algebras},  Represent. Theory  13  (2009), 349-370. 

\bibitem[FH]{FH}  W. Fulton, and J. Harris, {\it Representation Theory}, Graduate Texts in Mathematics, Springer, 2004.

\bibitem[GG]{GG} I. Gordon and  S. Griffeth, {\it Catalan numbers for complex reflection groups}, preprint arXiv:0912.1578v1 

\bibitem[GGOR]{GGOR}  V. Ginzburg, N. Guay, E. Opdam, and R. Rouquier, {\it On the category $\mathcal{O}$ for rational Cherednik algebra}, Invent. Math. Vol. 154, No. 3. (2003), 617-651.

\bibitem[H]{H} J. Humphreys, {\it Reflection groups and Coxeter groups}, Cambridge University Press, 1992.

\bibitem[M]{M} J. M\"uller. {\it Decomposition numbers for generic Iwahori-Hecke algebras of noncrystalographic type}, J. of Algebra Vol. 189 (1997), 125-149.

\bibitem[O]{O} E.M. Opdam, {\it A remark on the irreducible characters and fake degrees of finite real reflection groups}, Invent. Math. Vol. 120 (1995) 447-454. 

\bibitem[R1]{R1} R. Rouquier, {\it q-Schur algebras and complex reflection groups, I}, Moscow Math. J. Vol. 8, No 1 (2008), 119-158.

\bibitem[R2]{R2} R. Rouquier, {\it Representations of rational Cherednik algebras}, arXiv:math/0509252v2.

\bibitem[VV]{VV} Varagnolo, M. and E. Vasserot
{\it Finite-dimensional representations of DAHA and affine Springer fibers: The spherical case},
Duke Math. J. Vol. 147, No. 3 (2009), 439-540.
\end{normalsize}
\end{thebibliography}
\end{document}